\documentstyle[amstex]{article}
\numberwithin{equation}{section}
\newtheorem{thm}{Theorem}[section]
\newtheorem{lem}{Lemma}[section]
\newtheorem{defn}{Definition}[section]
\newtheorem{prop}{Proposition}[section]
\newtheorem{sublem}{Sublemma}[section]
\newtheorem{cor}{Corollary}[section]
\begin{document}
\Large
\title{ \LARGE \bf 
 Localized ASD moduli spaces based on the reduced cohomology group over Casson handles}
\author{ \LARGE \bf Tsuyoshi Kato}

\date{}
\maketitle

\begin{center}{ \bf  Introduction}\end{center}
One of the most important aspects in global analysis on non compact spaces 
is  the choice of functional spaces over them.
Of particular interest for us is analysis of elliptic differential operators 
over non compact manifolds, which  deeply touches topological structure  of their underlying spaces.
When one uses standard Sobolev spaces or weighted ones as the analytic setting,
 Fredholm property often breaks, which happens when  continuous spectrum
appears near zero.
Let us describe such situation in terms of the cohomology groups. Let $X$ be a space and:
\begin{gather*}
0 \;\longrightarrow \; C^0(X)
\;\stackrel{d_0}{\longrightarrow} \;
C^1(X)
  \stackrel{d_1}{\longrightarrow} \; 
C^2(X) \;
\longrightarrow \; 0.
\end{gather*}
be a bounded complex 
between  topological cochains $C^*(X)$.
One obtains  two different types of cohomology groups from this complex,
 where one is the ordinary cohomology $H^i(X)$ and the another
is the reduced one $\bar{H}^i(X)$:
$$H^i(X)= \text{ Ker } d_i  /  \text { im } d_{i-1}, \qquad
\bar{H}^i(X)= \text{ Ker } d_i  /  \bar{\text { im } d_{i-1}}$$
where $\bar{\text { im } d_{i-1}} $ is the closure of $\text { im } d_{i-1} \subset C^i(X)$.
There is  a canonical surjection
$H^i(X) \to  \bar{H}^i(X)$.

These coincide with each other when $X$ is compact,
while 
 let us consider 
the differential over $L^2$ functions on the universal covering space
of a compact manifold. It was verified by Brooks that 
these coincide with each other, if and only if the fundamental group is non amenable.
Such difference is reflected on the behavior of the spectrum of the Laplace operator near zero,
whether continuous spectrum attains zero or it is isolated.

In this paper we introduce a new functional analytic framework
which is based on the reduced cohomology groups of elliptic  complexes
over non compact manifolds, and 
develop a deformation theory of the Fredholm operators on  index theory.
Our original motivation to use 
 the  reduced cohomology arose  by the work by Ballmann, Br\"urning and Carron,
 which include index theory of elliptic differential operators 
 over cylindrical manifolds with boundary  ([BBC]).
 
We shall   apply it to non linear PDE analysis and obtain a  topological 
constrain of complexity of smooth structure over non compact
open subsets which are embedded into compact smooth four manifolds.

Let $E_i \to X$ be Euclidean vector bundles over $X$ and fix a sufficiently large $k$. 
Let us consider
a family of of  elliptic complexes:
\begin{gather*}
0 \;\longrightarrow \; L^2_{k+1}(E_0)  
\;\stackrel{d^0_{\mu}}{\longrightarrow} \; 
L^2_k(E_1)  
  \stackrel{d^1_{\mu}}{\longrightarrow} \; 
L^2_{k-1}(E_2) \;
\longrightarrow \; 0 \qquad (*)
\end{gather*}
  parametrized by   $ 0 \leq \mu \leq \delta$, and 
  introduce the Hilbert spaces   ${\frak L}_{k+1}(E_i)_{\mu}$ as the functional spaces
which we call the reduced Sobolev spaces with 
 their norms:
$$||u||^2_{{\frak L}_{k+1}(E_0)_{\mu}} =||d^0_{\mu}(u)||^2_{L^2_k}, \quad 
||w||^2_{{\frak L}_k(E_1)_{\mu}} = ||w'||^2_{L^2_k} + ||d^1_{\mu}(w'')||^2_{L^2_{k-1}}
$$
with respect to the orthogonal decomposition:
$$L^2_k(E_1) \cong V \oplus V^{\perp}, \quad w=w'+w''$$
where $V $ is the closure of the image $d^0_{\mu} (L^2_{k+1}(E_0)) \subset L^2_k(E_1)$.
We put
${\frak L}_{k-1} (E_2) = L^2_{k-1}(E_2)$ as the usual Sobolev space.

These functional spaces satisfy  some non standard properties, 
for example both $0$-th and $1$-st  cohomology groups are eliminated.
On the other hand
it turns out that both ${\frak L}_{k+1}(E_0)_{\mu}$ and ${\frak L}_k(E_1)_{\mu}$ 
admit the equivalent norms with
$L^2_{k+1}$ and $L^2_k$  respectively at $\mu$, 
 if it is of Fredholm such that  their cohomology groups
 $H^*_{\mu}$ satisfy $H^0_{\mu} =H^1_{\mu}=0$.
 
When the base space is non compact, 
the standard Sobolev spaces do not suffice to obtain Fredholm complexes
by elliptic PDE systems. Often we use smaller functional spaces by replacing them by 
the weighted  Sobolev spaces. A typical situation is the case of cylinderical manifolds.
 If we assign $\mu $ as the weight constants, then they constitute 
 the Fredholm complexes over cylinderical manifolds for small and positive $\mu >0$.
 In the case the functional spaces are the standard Sobolev spaces at $\mu =0$,
 where they are no more of Fredholm.
 
 Our first theorem is the following, which applies also to  the case of cylinderical manifolds.

\begin{thm} Suppose the above family of elliptic complexes $(*)$
are the  filtered Fredholm complexes 
for  $ 0 < \mu \leq \delta$
of non negative
 indices $m \geq 0$.
 Then  the family of the induced complexes  ${\frak C}_{\mu}$: 
\begin{gather*}
0 \;\longrightarrow \; {\frak L}_{k+1}(E_0)_{\mu}  
\;\stackrel{d^0_{\mu}}{\longrightarrow} \; 
{\frak L}_k(E_1)_{\mu}  
  \stackrel{d^1_{\mu}}{\longrightarrow} \; 
L^2_{k-1}(E_2) \;
\longrightarrow \; 0
\end{gather*}
are  also of Fredholm for all $0 \leq \mu \leq \delta$
whose  indices are  equal to:
$$ \dim H^2_{\mu} .$$

In particular if  $\dim H^0_{\mu} =0$ hold 
for all $0 < \mu \leq \delta$, 
then  the index of 
${\frak C}_0$
 is larger than or equal to
the Euler characteristics: 
$$\text{ ind } {\frak C}_0 \geq - \dim H^1_{\mu} + \dim H^2_{\mu}.$$
\end{thm}

In four manifold theory, 
Yang-Mills gauge theory is a fundamental tool to study topological  structure of smooth four manifolds.
It  uses the moduli space  which is given by the set of solutions to the ASD equation modulo gauge transformations.
In general, the infinitesimal structure of the ASD moduli space is based on  the 
Atiyah-Hitchin-Singer elliptic complex:
\begin{gather*}
0 \;\longrightarrow \; C^{\infty}(X)
\;\stackrel{d}{\longrightarrow} \;
C^{\infty}(X; \Lambda^1)
  \stackrel{d^+}{\longrightarrow} \; 
C^{\infty}(X; \Lambda^2_+) \;
\longrightarrow \; 0
\end{gather*}
The first differential 
$ C^{\infty}(X)
\;\stackrel{d}{\rightarrow} \;
C^{\infty}(X; \Lambda^1)$
corresponds to the infinitesimal gauge group action,
and closeness of the image corresponds to Hausdorff property in  the construction of the global
moduli space.
In the theory,   several  functional spaces have been used so far,
not only ordinary Sobolev spaces  based on  the ordinary cohomology theory,
but also  weighted Sobolev spaces.
There is another development by use of families of Banach spaces
in study of quasi-conformal mappings ([K4]).

In this paper we apply our functional analytic setting above, and 
introduce a new construction of the localized ASD moduli space
which is based on the reduced cohomology theory.
It turns out that this formalism makes the local construction of the moduli space quite canonical
in a situation when the differentials of AHS complex do not have closed range,
where the standard $L^2$ theory does not work directly.
We think  that our construction may have more chance to apply in another 
occasion, such as construction of instanton Floer homology groups ([F]).

Our main motivation is to study smooth complexity of Casson handles
inside smooth four manifolds.
Casson handle $CH(T)$ is an open four manifold with boundary,  which 
is homeomorphic to the standard open $2$ handle but far from diffeomorphic.
It is  parametrized by a signed and rooted infinite tree $T$, 
and 
if two such trees $T_1 \subset T_2$ admit embedding, 
then the corresponding Casson handles also admit smooth embedding 
$CH(T_2) \subset CH(T_1)$ in a reverse way,  preserving their attaching regions.
 So growth structure of the tree directly reflects complexity
of its smooth structure.

Casson handles
arise 
when a simply connected, oriented and smooth four manifold is
decomposed  topologically 
with respect to its intersection form (see [K3]).
Let $M$ be such a manifold with even type form.
 Then there exists an open subset $S \subset M$  homeomorphic to  the connected sums of 
$S^2 \times S^2$ removed $4$ cell, 
which is compatible with the decomposition of the form.
$S$  admits the induced  smooth structure and admits
 a smooth decomposition: 
$$S \cong D^4 \cup_{i=1}^{2l} CH(T_i).$$
So these  Casson handles are  embedded  inside 
$M$ smoothly, and 
 both  the end of $S$ and $S$ itself are simply connected.

We say that an  open four manifold $S$ has {\em tree-like end } if
there is a 
finite family of  signed trees $T_1, \dots, T_l$, such that
$S$ is diffeomorphic to $D^4 \cup^l_{j=1} CH(T_j)$, where
 every $CH(T_j)$ is attatched to the 
zero handle $D^4$ along 
the attatching $S^1$ of the first stage kinky handle which
corresponds to the root in $T_j$.

In [K1],  we have  introduced  a class of signed infinite trees which are
called {\em  trees of bounded type}.
Any tree  of bounded type grows polynomially.

Let $M$ be as above, and take an $SO(3)$ bundle $E \to M$.
A standard form  of the intersection form $ < \quad, \quad> = k(-E_8) \oplus lH$,
gives decomposition as
$H_2(M: {\mathbb Z}) \cong {\mathbb Z}^{8k} \oplus  {\mathbb Z}^{2l}$,
where $H$ is the hyperbolic $2$ by $2$ matrix.
We call such splitting  a marking of the form.
A marking gives an open four manifold $S = D^4 \cup^l_{j=1} CH(T_j) \subset M$. 
Notice that marking is not unique because of existence of lattice automorphisms.
The trees, and hence the Casson handles change if we choose different markings.
In [K3], we gave a proof of the following:
\begin{thm}
(1) For any marking on the K3 surface, the corresponding 
embedded Casson handles  cannot be all of bounded type.

(2) 
Let $M$ be K3 surface or  its logarithmic transforms $X_p$ for odd $p$.
Then there is an $SO(3)$ bundle $E \to M$ with $w_2(E) \ne 0$ 
which admits non empty generic markings,
so that
for any generic marking, the corresponding
embedded Casson handles  can not be all of bounded type.
\end{thm}
With  respect to the decomposition above, 
the second Stiefel-Whitney class $w_2$ of $E$  splits as
$w_2 =w_2^1 \oplus w_2^2$.
We  say that a marking is generic with respect to
the $SO(3)$ bundle $E$, 
if both $w_2^1 \ne 0$ and $w_2^2 \ne 0$ do not vanish.

Our aim is to construct another approach  to the proof of  theorem $0.2$ 
 by use of the new functional analytic setting described above.
 (1)   follows from (2)  with  [M] ([K3]).
 So we focus on the proof of (2) combining with [Kr] which concerns 
 existence of ASD connections over $M$.

  In [K1], we have explicitly  introduced complete Riemannian metrics of bounded geometry 
  on any Casson handles.

  \begin{thm}
 Let $S \cong D^4 \cup_{i=1}^{2l} CH(T_i)$ be the Riemannian-Casson handles
 of homogeneously bounded type.
  Then 
the  induced AHS complex:
\begin{gather*}
0 \;\longrightarrow \; {\frak L}_{k+1}(S)  
\;\stackrel{d}{\longrightarrow} \; 
{\frak L}_k(S;  \Lambda^1)  
  \stackrel{d^+}{\longrightarrow} \; 
L^2_{k-1}(S;  \Lambda^2_+) \;
\longrightarrow \; 0.
\end{gather*}
  is  of  Fredholm whose  index $I$ admits the bounds:
  $$l \leq I \leq 2l.$$
\end{thm}
This is a consequence of theorem $0.1$. 
In [K1], we have constructed weighted Sobolev spaces with the weights $0< \mu<<1$ 
such that the corresponding AHS complexes AHS$_{\mu}$ are filtered Fredholm. So
the family of the Fredholm complexes satisfies the assumption in theorem $0.1$.

Let us  describe the construction of our variant of 
Yang-Mills moduli theory.
Let $E \to S$ be an $SO(3)$ bundle whose trivialization near the end
is fixed.
A  connection $A$ over $E$ is called the {\em anti self dual}
(ASD),
if its curvature $F_A$ satisfies the equation:
$$F^+_A \equiv F_A +  * F_A=0.$$
The  construction of our localized ASD moduli space at $A$
uses the functional spaces ${\frak L}_k(A)$, however there causes  a problem 
if one tries to use the ASD equation itself, since it will not be well defined to formulate
the  self dual curvature form for elements
$A+a  \in {\frak L}_k(A)$. So we use a kind of regularization of the self dual curvature
which replaces $F^+$ by $\tilde{F}^+$,  by use of the spectral decomposition.
Let us cut off the spectra near zero 
as $P_{f_{\epsilon}} = \int_0^{\infty} f_{\epsilon}(\lambda) d E(\lambda)$,
where $f_{\epsilon}$ vanishes on $ [0,\epsilon] $.
Then the regularization is given by: 
$$ \tilde{F}^+(A+a)  =  d_A^+(a) +  
(Q_{f_{\epsilon}}(a) \wedge Q_{f_{\epsilon}}(a))^+$$
where $Q_{f_{\epsilon}}(a) = (d^+_A)^* \Delta_A^{-1} P_{f_{\epsilon}}( d^+_A(a))$.
Then we define  the localized ASD moduli space by:
 $${\frak M}_k(A) =\{ A + a :  \ \tilde{F}^+(A+a)=0 , \ \  a \in (\text{ Ker } d^+_A)^{\perp} \subset {\frak L}_k(A)\}.$$

 Let us describe the  idea of our approach to (2) in   theorem $0.2$ roughly.
 Let $M$ be a simply connected,  oriented,  closed  and smooth four manifold of even type,
 equipped with a marking:
$$ \Phi: (H_2(M; {\mathbb Z}), <,>)  \cong
(\oplus^{8k} {\mathbb Z}  \oplus^{2l} {\mathbb Z},  \quad k(-E_8)  \oplus lH)$$
where $k \geq 2$ and $l \geq 3$.
By Casson-Freedman theory, one finds an open four manifold $S$ 
with tree-like end, homeomorphic to the interior of  
$l(S^2 \times S^2) \backslash D^4$,
and finds  a smooth embedding
$S \hookrightarrow M$ which induces an embedding of the
form:
 $$(H^2(S:{\mathbb Z}), < \quad, \quad >) \cong (\oplus^{2l} {\mathbb Z}, lH )
  \hookrightarrow (H^2(M:{\mathbb Z}), <\quad, \quad>).$$

Let us assume that  the Donaldson's invariant is  non zero, and suppose
   the embedded Casson handles could be  all of bounded type.
    Let us  equip with the Riemannian metric $g$ on $S$ 
in  theorem $0.3$.
 We induce a  contradiction as below.
  Let us choose an exaustion on $S$ by compact subsets as
 $K_0 \subset K_1 \subset \dots \subset S \subset M$.
 Choose a family of generic Riemannian metrics $\{g_i\}_i$ on $M$
 such that $g_i|K_i \sim g|K_i$ are sufficiently near each other in
$C^{\infty}$.
Take a family of  ASD connections $A_i$ with respect to $(M,g_i)$
 which  converges to an $L^2$
 ASD connection $A$ with respect to $(S,g)$.
  Then we obtain the non empty ASD moduli space
 $A \in {\frak M}_k(A)$
 over $S$, where we use the new functional spaces.
$A$ may  not be  regular, and so we use 
perturbation of the regularized ASD equation, 
which is the simplest kind where we do not care about gauge group actions.
 One can find a solution to the perturbed equation near $A$
in  the functional space ${\frak L}_k(A)$.
  The index computation 
 shows that the dimension should be negative,
 which gives a contradiction.

M.Tsukamoto pointed out to me that $L^2$ ASD connections can not be trivialized 
near infinity over  general open four manifolds, and 
our argument in proposition $1.1$ [K3] is not enough,
 to whom the author thanks.
We also  fill in  the detail of this gap in section $3$, 
and verify that the argument  still  works in the case of our Riemannian-Casson handles.
We shall introduce a notion of $T^* \times {\mathbb N}$ structure on Riemannian manifolds,
and verify that Casson handles admit such structure, and those Riemannian manifolds with
a condition on the fundamental group admit such trivialization.

 Therorem $0.2$ implies that 
theose Casson handles in K3 should grow much more than 
bounded type so that our framework of the Fredholm theory also
breaks. On the other hand
 the $L^2$ ASD connection $A$ 
 over the Riemannian-Casson  handles inside $K3$ surface certainly exist by the 
 above deformation process of the Riemannian metrics.
It follows from our argument that 
 the cokernel of:
 $$d^+_A: {\frak L}_k(A) \to L^2_{k-1}(S: Ad(P) \otimes \Lambda^2_+)$$
should be $0$ or infinite dimensional.

In  the former case, the ASD moduli space
near $A$ would consist of zero dimensional regular smooth manifold.
However since the growth of the Casson handle inside $K3$ surface
would be so wild, it seems reasonable to predict its behavior as:
\vspace{3mm} \\
{\em Conjecture 0.1:}
Let $A$ be the $L^2$ ASD connection over the Riemannian-Casson  handle inside $K3$ surface.
Then the cockerel of
 $d^+_A$ will be of  infinite dimension.
\vspace{3mm}

From geometric analysis view point, it would be of interest  for us
to develop construction of the  global ASD moduli space and the
 gauge group action over our functional  spaces.

\section{Reduced cohomology group}
{\bf 1.A Functional spaces:}
Let $(X,g)$ be a complete Riemannian manifold, and
$E_i \to X$ be $SO(N)$ vector bundles
which are trivialized near infinity  for $i =1,2,3$.

Suppose there is a complex between the Sobolev spaces:
\begin{gather*}
0 \;\longrightarrow \; L^2_{k+1}(E_0)  
\;\stackrel{d^0}{\longrightarrow} \; 
L^2_k(E_1)  
  \stackrel{d^1}{\longrightarrow} \; 
L^2_{k-1}(E_2) \;
\longrightarrow \; 0.
\end{gather*}
by elliptic differential operators.
\vspace{3mm} \\
{\em Example 1.1:}
Our main application is the
Atiyah-Hitchin-Singer complex:
\begin{gather*}
0 \;\longrightarrow \; L^2_{k+1}(X)  
\;\stackrel{d}{\longrightarrow} \; 
L^2_k(X; \Lambda^1)  
  \stackrel{d^+}{\longrightarrow} \; 
L^2_{k-1}(X; \Lambda^2_+) \;
\longrightarrow \; 0
\end{gather*}
over  smooth four  manifolds, 
where $d^+$ is the composition of $d$ with the projection to the self dual part on $2$ form.

\vspace{3mm}

Let us  study analytic behavior of deformation of these elliptic complexes.
Let $\delta >0$ be a small and positive number, 
and consider a smooth deformation of the elliptic differential operators:
\begin{gather*}
0 \;\longrightarrow \; L^2_{k+1}(E_0)  
\;\stackrel{d_{\mu}^0}{\longrightarrow} \; 
L^2_k(E_1)  
  \stackrel{d^1_{\mu}}{\longrightarrow} \; 
L^2_{k-1}(E_2) \;
\longrightarrow \; 0.
\end{gather*}
 for  $ 0 \leq \mu \leq \delta$.
 Let us denote  by $H^*_{\mu}$ as
 their (unreduced) cohomology groups.

Notice that the differential:
$$d^0: L^2_{k+1}(E_0) \to L^2_k(E_1).$$
 does not have closed range in general.
 For example 
  the differentials on functions do not have closed range over 
 the cylinderical manifolds or $ {\mathbb R}^n$.

Let us  introduce  new functional spaces which we call the reduced Sobolev spaces.
\begin{defn}
(1)
The reduced Sobolev space ${\frak L}_{k+1}(E_0)$
(which depend on $\mu$) is
given 
by the maximal extension of the domain
$C^{\infty}_c(E_0)$ of $d_{\mu}^0$ 
 with  the norm:
$$||u||^2_{{\frak L}_{k+1}} =||d_{\mu}^0(u)||^2_{L^2_k}.$$

(2)
 ${\frak L}_k(E_1)$ is
 given 
by the closure of $L^2_k(E_1)$ 
 with  the norm:
$$||w||^2_{{\frak L}_k} = ||w'||^2_{L^2_k} + ||d^1_{\mu}(w'')||^2_{L^2_{k-1}}
$$
with respect to the orthogonal decomposition $w=w'+w''$ as:
$$L^2_k(E_1) = {\frak L}_{k+1}(E_0) \oplus {\frak L}_{k+1}(E_0)^{\perp}.$$ 
\end{defn}
Notice that ${\frak L}_{k+1}(E_0) $ is identified with the closed linear subspace 
 $d^0_{\mu}({\frak L}_{k+1}(E_0)) \subset  L^2_k(E_1)$.
Our choice of the norms induces the 
 family of  bounded complexes  ${\frak C}_{\mu}$:
\begin{gather*}
0 \;\longrightarrow \; {\frak L}_{k+1}(E_0)  
\;\stackrel{d^0_{\mu}}{\longrightarrow} \; 
{\frak L}_k(E_1)  
  \stackrel{d^1_{\mu}}{\longrightarrow} \; 
L^2_{k-1}(E_2) \;
\longrightarrow \; 0.
\end{gather*}

\begin{lem}
Suppose $X$ is compact without boundary.
Then 
 there are canonical isomorphisms for $*=0,1$:
$${\frak L}_l(E_*) \cong L^2_l(E_*) / H^*_{\mu}$$

In particular
there are embeddings:
$${\frak L}_l(E_*) \hookrightarrow L^2_l(E_*).$$
\end{lem}
{\em Proof:}
Let us consider the case $*=0$. 
There is a continuous map from
$L^2_{k+1}(E_0)$ to ${\frak L}_{k+1}(E_0)$, since 
 a priori estimates hold:
$$||v||_{L^2_{k+1}} \geq C ||v||_{{\frak L}_{k+1}}.$$

By the assumption, the spectrum of $\Delta^0 \equiv (d^0_{\mu})^* d^0_{\mu}$ is discrete, and 
let us decompose $L^2_{k+1}(E_0) = {\mathbb R}^m \oplus V$, 
where $V = $ Ker $(d^0_{\mu})^{\perp}$.
Then the estimates:
$$||v||_{L^2_{k+1}} \leq C' ||d^0_{\mu}(v)||_{L^2_k} = C' ||v||_{{\frak L}_{k+1}}$$
 hold  for any $v \in V$ and
  for some constant $C'$.
So the  continuos map 
 from  $L^2_{k+1}(X)$ to ${\frak L}_{k+1}$ is surjective and hence has closed range
with kernel $= {\mathbb R}^m$.

Next consider $*=1$ case.
By definition of the norm, 
${\frak L}_k(E_1) $ splits as the direct sum of 
$d^0_{\mu}({\frak L}_{k+1}(E_0) ) \subset L^2_k(E_1)$ 
with $d^1_{\mu}({\frak L}_k(E_1)) \subset L^2_{k-1}(E_2)$,
both of which 
consist of the closure of the images of $d^0_{\mu}(L^2_{k+1}(E_0))$ and 
$d^1_{\mu}(L^2_k(E_1))$ respectively.
By the assumption, both $d^0_{\mu}$ and $d^1_{\mu}$ have closed range, and hence
it splits as: 
$$ {\frak L}_k(E_1) \cong d^0_{\mu}(L^2_{k+1}(E_0) ) \oplus  d^1_{\mu}(L^2_k(E_1)) .$$
So it is enough to see the isomorphism:
$$L^2_k(E_1) / H^1_{\mu} \cong d^0_{\mu}(L^2_{k+1}(E_0) ) \oplus  d^1_{\mu}(L^2_k(E_1)) $$
where we identify $H^1_{\mu} \cong {\mathbb R}^m \subset L^2_k(E_1)$.

Let $V = (  d^0_{\mu}(L^2_{k+1}(E_0) ) \oplus {\mathbb R}^m)^{\perp} \subset L^2_k(E_1)$.
Then as in the case of $*=0$, $d^1_{\mu}: V \cong d^1_{\mu}(L^2_k(E_1))$ gives the isomorphism.
So we obtain the isomorphism:
$$L^2_k(E_1) \cong  d^0_{\mu}(L^2_{k+1}(E_0) ) \oplus   d^1_{\mu}(L^2_k(E_1)) \oplus  {\mathbb R}^m$$
which gives the desired isomorphism.

This completes the proof.

\vspace{3mm}

Below let us state an application of lemma $1.1$.
Let us say that an elliptic complex is reversible, if for any compact subset $K  \subset X$, 
there is another compact submanifold $K \subset K'  \subset X$
such that the restriction of the  complex over $K'$ can be extended to
another elliptic complex over the double $DK'$.

\begin{cor} Suppose $X$ is possibly non compact,
and the elliptic complex is reversible. 
Let $\varphi: X \to [0,1]$ be a cut off function with compact support.

Then there admit continuous functional:
$$\varphi  : {\frak L}_l(E_*)  \to {\frak L}_l(E_*)_c$$
which are induced from  multiplication by $\varphi$.
\end{cor}
{\em Proof:}
Let $K \subset X$ be a compact subset 
which contains support of $\varphi$ in its interior.
Let  $  {\frak L}_l(E_*)_K \subset   {\frak L}_l(E_*)$ be the closure
of the image of $L^2_l(E_*|K)_0$, and $DK $ be the double of $K$, which is a closed manifold.
One can extend the differentials so that anoter elliptic complexes:
\begin{gather*}
0 \;\longrightarrow \; L^2_{k+1}(E_0|DK)  
\;\stackrel{(d_{\mu}^0)'}{\longrightarrow} \; 
L^2_k(E_1|DK)  
  \stackrel{(d^1_{\mu})'}{\longrightarrow} \; 
L^2_{k-1}(E_2|DK) \;
\longrightarrow \; 0.
\end{gather*}
are obtained,
 which extend the restriction of the original ones over $K$.
Then we obtain the embeddings by lemma $1.1$:
$$\Phi: {\frak L}_l(E_*)_K  \hookrightarrow {\frak L}_l(E_*|DK) \hookrightarrow L^2_l(E_*|DK)$$
and the composition of the following gives the desired map:
$$\varphi(a) \equiv \varphi \Phi(a) \in L^2_l(E_*)_K \subset L^2_l(E_*)_c \to {\frak L}_k(E_*)_c.$$
 
 This completes the proof.

\vspace{3mm}

Let $X$ be a complete Riemannian manifold.

\begin{lem}
(1)
 There are isometries of the Hilbert spaces:
 \begin{align*}
& d_{\mu}^0 : {\frak L}_{k+1} (E_0) \cong  d_{\mu}^0({\frak L}_{k+1}(E_0)) 
\subset L^2_k(E_1) . \\
& {\frak L}(E_1) \cong d_{\mu}^0({\frak L}_{k+1}(E_0)) \oplus d^1_{\mu}({\frak L}(E_1))
\subset L^2_k(E_1) \oplus L^2_{k-1}(E_2).
\end{align*}

(2)
Suppose: 
\begin{gather*}
0 \;\longrightarrow \; L^2_{k+1}(E_0)  
\;\stackrel{d_{\mu}^0}{\longrightarrow} \; 
L^2_k(E_1)  
  \stackrel{d_{\mu}^1}{\longrightarrow} \; 
L^2_{k-1}(E_2) \;
\longrightarrow \; 0.
\end{gather*}
is a Fredholm complex with $H^0_{\mu}=H^1_{\mu}=0$ at $\mu$.

Then there are constants $C_{\mu}  >0$ such that the estimates hold:
\begin{align*}
& C^{-1}_{\mu} ||u||^2_{L^2_{k+1}} \leq ||u||^2_{{\frak L}_{k+1}(E_0)} \leq   C_{\mu} ||u||^2_{L^2_{k+1}}, \\
&  C^{-1}_{\mu} ||w||^2_{L^2_k} \leq ||w||^2_{{\frak L}_k(E_1)} \leq   C_{\mu} ||w||^2_{L^2_k}.
\end{align*}
\end{lem}
{\em Proof:}
(1) follows by definition (see the proof of lemma $1.1$).

For (2), let us verify the first equivalence.
The uniform bounds:
$$C^{-1}_{\mu}  ||u||_{L^2_{k+1}} \leq 
 ||d_{\mu}^0(u)||_{L^2_k} \leq  C_{\mu}||u||_{L^2_{k+1}} $$
hold by the assumption, and so 
the conclusion holds.

Next let us consider the second one.
If $w= d_{\mu}^0(u)$ for some $u \in {\frak L}_{k+1}(E_0)$, then
the estimates hold by definition of the norm.

Let $w \in (d^0_{\mu}({\frak L}_{k+1}(E_0)))^{\perp}$.
Then the estimates:
$$C_{\mu}^{-1} ||d^1_{\mu}(w)||_{L^2_{k-1}}
\leq ||w||_{L^2_k} \leq C_{\mu}||d^1_{\mu}(w)||_{L^2_{k-1}}$$ holds
since  $H^1_{\mu}(X)=0$ and Fredholmness of the complex.
Any element $w$ in ${\frak L}_k(E_1)$ can be 
given by the direct sum  of these cases. This  verifies the second inequalities.
This completes the proof.
\vspace{3mm} \\
{\bf 1.A.2 Filtration of functional spaces:}
Let us consider a smooth family   of elliptic complexes
 for $0 \leq \mu \leq \delta$:
\begin{gather*}
0 \;\longrightarrow \; L^2_{k+1}(E_0)  
\;\stackrel{d_{\mu}^0}{\longrightarrow} \; 
L^2_k(E_1)  
  \stackrel{d^1_{\mu}}{\longrightarrow} \; 
L^2_{k-1}(E_2) \;
\longrightarrow \; 0.
\end{gather*}

\begin{defn}
Let us say that the family of complexes are {\em filtered}, if
there are infinite embeddings by Hilbert spaces:
$$C^{\infty}_c(E_*) \subset W_{\mu}(E_*) \subset W_{\mu'}(E_*) \subset 
W_0 = L^2_l(E_*)$$
for any $0 \leq \mu' \leq \mu \leq \delta$ such that
there are isomorphisms:
$$I_{\mu} : L^2_l(E_*) \cong W_l(E_*)_{\mu}$$
with $I_0 =id$, 
which are compatible with their complexes:
\begin{gather*}
0 \;\longrightarrow \; W_{k+1}(E_0)_{\mu}  
\;\stackrel{d^0}{\longrightarrow} \; 
W_k(E_1)_{\mu}  
  \stackrel{d^1}{\longrightarrow} \; 
W_{k-1}(E_2)_{\mu} \;
\longrightarrow \; 0
\end{gather*}
where both $d^0$ and $d^1$ are independent of $\mu$.
\end{defn}
We will describe in $1.B$
that AHS complexes with respect to the weighted Sobolev spaces are filtered
over cylindrical manifolds.

\begin{prop}
Let us consider the filtered  complexes  for $0 \leq \mu \leq \delta$,
which are of Fredholm
of the same indices 
with cohomology groups $H^0_{\mu}=
H^1_{\mu}=0$ and $\dim H^2_{\mu}= m$
for all $0 < \mu \leq \delta$.

Let us consider the induced  complexes ${\frak C}_{\mu}$:
\begin{gather*}
0 \;\longrightarrow \; {\frak L}_{k+1}(E_0)_{\mu}  
\;\stackrel{d_{\mu}^0}{\longrightarrow} \; 
{\frak L}_k(E_1)_{\mu}  
  \stackrel{d^1_{\mu}}{\longrightarrow} \; 
L^2_{k-1}(E_2)  \;
\longrightarrow \; 0.
\end{gather*}
Then it is of Fredholm whose cohomolgy groups
${\bf H}^i_{\mu}$  satisfy: 
 $${\bf H}^0_{\mu} = {\bf H}^1_{\mu}=0, \quad 
 \dim {\bf H }^2_{\mu}=m$$
 for all $0 \leq \mu \leq \delta$.
 In particular it is also Fredholm at $\mu =0$.
 
 Moreover there is $m$ dimensional linear subspace 
 $V \subset L^2_{k-1}(E_2)$ such that
  the projections:
$$\pi : d^1_{\mu}({\frak L}_k(E_1)_{\mu})  \to L^2_{k-1}(E_2) / V$$ are onto
for all small $0 \leq \mu $.
\end{prop}
{\em Proof:} 
{\bf Step 1:}
By definition,  ${\bf H}^0_{\mu}=0$ and ${\bf H}^1_{\mu}=0$ hold
for all $0 \leq \mu \leq \delta$.

{\bf Step 2:}
Suppose $H^2_{\mu}=0$ for all positive $ \mu >0$.
For  any $v \in L^2_{k-1}(E_2)$,
there is  $w \in L^2_k(E_1)$ with $d^1_{\mu}(w)=v$.
One can regard  $w \in {\frak L}_k(E_1)_{\mu}$ by definition of the norm,
and so ${\bf H}^2_{\mu}=0$ hold for all $0 < \mu \leq \delta$.

Let us  verify 
${\bf H}^2_0 =0$.
Notice that  
$d^1_{\mu} : {\frak L}^2_k(E_1)_{\mu}  
\to L^2_{k-1}(E_2) $
has closed range for all $0 \leq  \mu \leq \delta$
by definition of the norm.
So
it is enough to see
that the image  is dense.
Let:
\begin{gather*}
0 \;\longrightarrow \; W_{k+1}(E_0)_{\mu}  
\;\stackrel{d^0}{\longrightarrow} \; 
W_k(E_1)_{\mu}  
  \stackrel{d^1}{\longrightarrow} \; 
W_{k-1}(E_2)_{\mu} \;
\longrightarrow \; 0
\end{gather*}
be the filtered complexes.
Let us take any $v \in L^2_{k-1}(E_2)$,
and choose  approximations 
$v_{\mu} \in W_{k-1}(E_2)_{\mu}$
which converge to $v $ in $L^2_{k-1}(E_2)$.
Since ${\bf H}^2_{\mu} =0$ for all positive $\mu >0$,
there are $w_{\mu} \in W_k(E_1)_{\mu}  $
with $d^1(w_{\mu} )=v_{\mu}$.
In particular the image is dense as desired.

{\bf Step 3:}
Since  the Fredholm indices are invariant 
under continuous deformations,
$H^2_{\mu}$ have constant rank $m$ for all positive $ \mu >0$.

By definition,  $d^1_{\mu}({\frak L}_k(E_1)_{\mu})  \subset L^2_{k-1}(E_2)$
are closed subspaces for all $0 \leq \mu \leq \delta$.
We claim that there is a vector subspace $V \subset L^2_{k-1}(E_2)$
of dimension $m$
such that the projections:
$$\pi : d^1_{\mu}({\frak L}_k(E_1)_{\mu})  \to L^2_{k-1}(E_2) / V$$ are onto
and hence isomorphic
for all positive $\mu >0$.
Notice that  the orthogonal complement $V_{\mu} = d^1_{\mu}({\frak L}_k)^{\perp}$
 in $L^2_{k-1}(E_2)$
is a    smooth  family of $m$ dimensional vector subspaces.

Suppose contrary. Then for any $m$ dimensional vector space $V$, 
there is a small $\mu_0 >0$ such that
$\pi$ are onto for $\mu > \mu_0$ but not the case at $\mu_0$.
This implies that there is a line $l \subset V$ such that $l$ is contained in the image
of $d^1_{\mu_0}$.
So for any $\mu >0$, there is smaller $\mu > \mu_0$ so that
$V_{\mu}$ contain a line in $V_{\mu_0}^{\perp}$, which cannot happen.

Notice that 
the closure of $d^1( L^2_k(E_1)) $ is equal to
$ d^1({\frak L}_k(E_1))$ at $\mu=0$.
Let $W \subset L^2_{k-1}(E_2)$ be the co-kernel 
of $ d^1({\frak L}_k(E_1))$.
We claim that dimension of $W$ does not exceed $m$.
Let us consider the extension of the complexes:
\begin{gather*}
0 \;\longrightarrow \; L^2_{k+1}(E_0)  
\;\stackrel{(d_{\mu}, 0) }{\longrightarrow} \; 
L^2_k(E_1)  \oplus V
  \stackrel{d^+_{\mu} + \text{ id } }{\longrightarrow} \; 
L^2_{k-1}(E_2) \;
\longrightarrow \; 0.
\end{gather*}
This is acyclic for $0 < \mu \leq \delta$, 
and one can obtain the filtered complexes:
\begin{gather*}
0 \;\longrightarrow \; W_{k+1}(E_0)_{\mu}  
\;\stackrel{(d^0,0) }{\longrightarrow} \; 
W_k(E_1)_{\mu}   \oplus I_{\mu}(V)
  \stackrel{d^1 + \text{ id }}{\longrightarrow} \; 
W_{k-1}(E_2)_{\mu} \;
\longrightarrow \; 0
\end{gather*}
Let us take any $v \in L^2_{k-1}(E_2)$.
Then for any $\epsilon >0$, there is some $\mu >0$ and $v_{\mu} \in 
W_{k-1}(E_2)_{\mu}$
so that $||v - v_{\mu} ||_{L^2_{k-1}} < \epsilon$
with $v_{\mu} \in  d^1(W_k(E_1)_{\mu}) + I_{\mu}(V)$,
where dimension of $I_{\mu}(V)=m$.
This verifies the claim.

{\bf Step 4:}
The continuous family of the complexes ${\frak C}_{\mu}$:
\begin{gather*}
0 \;\longrightarrow \; {\frak L}^2_{k+1}(E_0)_{\mu}  
\;\stackrel{d^0_{\mu}}{\longrightarrow} \; 
{\frak L}^2_k(E_1)_{\mu}  
  \stackrel{d^1_{\mu}}{\longrightarrow} \; 
L^2_{k-1}(E_2) \;
\longrightarrow \; 0.
\end{gather*}
are Fredholm for all $0 \leq \mu \leq \delta$ by step $3$.

The following abstract  lemma  finishes the proof of proposition $1.1$:
\begin{lem} 
Let $V_t \subset H$ be a smooth  family of cofinite dimensional
vector subspaces. Suppose 
$$D_t: V_t \to H^2$$
is a family of uniformly  bounded maps 
with closed range for $t \in [0,1]$
such that they are isomorphic for all $t \in (0,1]$.
If $D_0$ is surjective, 
 then $D_0$ is injective.
\end{lem}
{\em Proof:}
Suppose contrary, and consider the isomorphism  
$D_0:  V_0 /  \text{ ker } D_0 \cong H_2$, and extend it as:
$$D_t:  V_t /   P_t(\text{ ker } D_0) \cong H_2$$ for small $t \in [0, \epsilon )$,
where $P_t$ are the projections to $V_t$.
It should be a family of isomorphisms,
since it is an open condition, which cannot happen.

This completes the proof.

\vspace{3mm}

Now we verify the following:
\begin{thm}
Let us consider the   filtered  complexes for $0 \leq  \mu \leq \delta$
which are of Fredholm of non negative
 indices $m \geq 0$ for $ 0 < \mu \leq \delta$.
 
 Then  the family of the induced complexes  ${\frak C}_{\mu}$: 
\begin{gather*}
0 \;\longrightarrow \; {\frak L}^2_{k+1}(E_0)_{\mu}  
\;\stackrel{d^0_{\mu}}{\longrightarrow} \; 
{\frak L}^2_k(E_1)_{\mu}  
  \stackrel{d^1_{\mu}}{\longrightarrow} \; 
L^2_{k-1}(E_2) \;
\longrightarrow \; 0
\end{gather*}
are  also of Fredholm for all $0 \leq \mu \leq \delta$
whose  indices are  equal to
$\dim H^2_{\mu} $.
\end{thm}
{\em Proof:}
{\bf Step 1:}
We have seen the conclusion in proposition $1.1$
for the special case when 
 $H^0_{\mu}=
H^1_{\mu}=0$ hold
for all $0 < \mu \leq \delta$.

Let us 
 consider the Fredholm complexes for $\mu >0$:
\begin{gather*}
0 \;\longrightarrow \; L^2_{k+1}(E_0)  
\;\stackrel{d_{\mu}}{\longrightarrow} \; 
L^2_k(E_1)  
  \stackrel{d^+_{\mu}}{\longrightarrow} \; 
L^2_{k-1}(E_2) \;
\longrightarrow \; 0.
\end{gather*}
Let us choose finite dimensional vector spaces:
$$V_{\mu}^0 = \text{ Ker } d_{\mu} \subset L^2_{k+1}(E_0), \quad
V_{\mu}^1  \subset \text{ Ker } d^1_{\mu} \subset L^2_k(E_1), \quad
V_{\mu}^2 \subset L^2_{k-1}(E_2)$$ 
which represent $H^0_{\mu}, H^1_{\mu}$ and $H^2_{\mu}$ respectively.

{\bf Step 2:} For $\mu >0$, 
suppose dim $V_{\mu}^2 $ is less than or equal to the Fredholm index.
 Then $\dim V_{\mu}^0 \geq \dim V_{\mu}^1$ holds.
 We verify that the index of ${\frak C}_{\mu}$  is 
  equal to $\dim H^2_{\mu}$.

 Let us prepare another vector space $W$ with the acyclic complexes between finite dimensional spaces:
 \begin{gather*}
0 \;\longrightarrow \; V_{\mu}^0 
\;\stackrel{f_{\mu}}{\cong} \; 
V_{\mu}^1 \oplus W
  \stackrel{}{\longrightarrow} \; 
0 \;
\longrightarrow \; 0.
\end{gather*}
Then one can add the extra vector spaces in the complex:
 \begin{align*}
0 \;\longrightarrow \; L^2_{k+1}(E_0)   = (V_{\mu}^0)^{\perp} \oplus V_{\mu}^0 &
\;\stackrel{d_{\mu} \oplus f_{\mu}}{\longrightarrow} \; 
L^2_k(E_1)  \oplus W=  (V_{\mu}^1)^{\perp} \oplus V_{\mu}^1 \oplus W \\
&   \stackrel{d^+_{\mu}}{\longrightarrow} \; 
L^2_{k-1}(E_2)  \;
\longrightarrow \; 0.
\end{align*}

Let 
$I_{\mu} : L^2_l(E_*) \cong W_l(E_*)_{\mu}$
be the isomorphisms, and put
$\tilde{f}_{\mu} =  I_{\mu}^{-1} f_{\mu} I_{\mu}$.
Then consider:
\begin{gather*}
0 \;\longrightarrow \; W_{k+1}(E_0)_{\mu}  
\;\stackrel{d^0 \oplus \tilde{f}_{\mu}}{\longrightarrow} \; 
W_k(E_1)_{\mu}  \oplus W
  \stackrel{d^1}{\longrightarrow} \; 
W_{k-1}(E_2)_{\mu}   \;
\longrightarrow \; 0
\end{gather*}
which is  Fredholm  complex
with $H^0_{\mu} =0$ and $H^1_{\mu}=0$, and
the index is equal to $\dim H^2_{\mu}$.

Now the induced complex:
$$
0 \;\longrightarrow \; {\frak L}_{k+1}(E_0)_{\mu}    \oplus V_{\mu}^0 
\;\stackrel{d_{\mu} \oplus f_{\mu}}{\longrightarrow} \; 
{\frak L}_k(E_1)_{\mu}  \oplus  V_{\mu}^1 \oplus W 
   \stackrel{d^+_{\mu}}{\longrightarrow} \; 
L^2_{k-1}(E_2)  \;
\longrightarrow \; 0.
$$
is  of Fredholm,  
which is chain homotopy equivalent to the original induced complex.
Its  index is equal to  $\dim H^2_{\mu}$ by lemma $1.2$.

{\bf Step 3:}
Suppose $\dim V_{\mu}^2 $ is larger than  the Fredholm index.
Then
the inequality  $\dim V_{\mu}^0 \leq  \dim V_{\mu}^1$ holds.
 Let us take another vector space $W$ with the acyclic complex between finite dimensional spaces:
 \begin{gather*}
0 \;\longrightarrow \; V_{\mu}^0 
\;\stackrel{f_{\mu}}{\longrightarrow} \; 
V_{\mu}^1 
  \stackrel{g_{\mu}}{\longrightarrow} \; 
W \;
\longrightarrow \; 0.
\end{gather*}
Then
 one can add the extra vector spaces in the complex:
 \begin{align*}
0 \;\longrightarrow \; L^2_{k+1}(E_0)   = (V_{\mu}^0)^{\perp} \oplus V_{\mu}^0 &
\;\stackrel{d^0_{\mu} \oplus f_{\mu}}{\longrightarrow} \; 
L^2_k(E_1) =  (V_{\mu}^1)^{\perp} \oplus V_{\mu}^1 \\
&   \stackrel{d^1_{\mu}\oplus g_{\mu} }{\longrightarrow} \; 
L^2_{k-1}(E_2) \oplus W \;
\longrightarrow \; 0.
\end{align*}
The index of this complex is 
$ \dim H^2_{\mu} $.
By the same way as step $2$, 
 the induced complex:
\begin{gather*}
0 \;\longrightarrow \; {\frak L}^2_{k+1}(E_0)_{\mu}  \oplus V^0_{\mu}
\;\stackrel{d^0_{\mu} \oplus f_{\mu} }{\longrightarrow} \; 
{\frak L}^2_k(E_1)_{\mu} \oplus V^1_{\mu}
  \stackrel{d^1_{\mu} \oplus g_{\mu}}{\longrightarrow} \; 
L^2_{k-1}(E_2) \oplus W\;
\longrightarrow \; 0.
\end{gather*}
is of Fredholm  of the  index $\dim H^2_{\mu}$,
which is chain homotopy equivalent to the original induced complex.

{\bf Step 4:}
Let us verify that $ \dim H^2_{\mu}$ are constant for all small $0< \mu$.
In fact there is a finite dimensional vector subspace 
$V \subset W_{k-1}(E_2)_{\mu} \subset W_{k-1}(E_2)_{\mu'}$
such that for any $v  \in W_{k-1}(E_2)_{\mu}$, there is some $a \in V$ such that
$v-a = d^1(w)$ holds for some $w \in W_k(E_1)_{\mu}$.

So  $ \dim H^2_{\mu} \geq  \dim H^2_{\mu'}$ hold for any $\mu \geq \mu'$, since 
 any elements in  $W_{k-1}(E_2)_{\mu'}$ can be approximated by 
elements in  $W_{k-1}(E_2)_{\mu}$.

{\bf Step 5:}
Combining with step $2,3,4$, it follows that the indices of 
the induced  complexes are  equal to $ \dim H^2_{\mu}$ for all
$0 < \mu \leq \delta$.
It follows from 
step $3$ and step $4$  in proposition $1.1$ that 
 it is a  family of Fredholm complexes for  all $0 \leq   \mu \leq \delta$,
whose indices coincide with 
$ \dim H^2_{\mu} $.

This completes the proof.
\vspace{3mm}

\begin{cor}
Consider the situation in theorem $1.1$.

If $\dim H^0_{\mu} =0$ hold 
for all $0 < \mu \leq \delta$, 
then  ${\frak C}_{\mu}$
are   of Fredholm for all $0 \leq \mu \leq \delta$
whose  indices are larger than or equal to 
the Euler characteristics: 
$$ - \dim H^1_{\mu} + \dim H^2_{\mu}
\leq \text{ ind } {\frak C}_{\mu} = \dim H^2_{\mu}.$$
\end{cor}
$\quad$ 
\vspace{3mm} \\
{\bf 1.B AHS complexes over cylindrical manifolds:}
The Atiyah-Hitchin-Singer complex is the elliptic differential complex over a
Riemannian  four manifold $X$:
\begin{gather*}
0 \;\longrightarrow \; L^2_{k+1}(X)  
\;\stackrel{d}{\longrightarrow} \; 
L^2_k(X; \Lambda^1)  
  \stackrel{d^+}{\longrightarrow} \; 
L^2_{k-1}(X; \Lambda^2_+) \;
\longrightarrow \; 0
\end{gather*}
where $d^+$ is the composition of the differential with the projection
to the self dual $2$ forms.

Let $X$ be a complete Riemannian manifold such that 
it is isometric to the product $M \times [0, \infty)$ except a compact subset $K \subset X$.
Such  space is called as a cylindrical manifold.

In $1.B$ we verify  the following:
\begin{prop} 
Let $X$ be a cylindrical four manifold.

(1) There is a filtered AHS complexes over $X$ and  positive 
$\delta >0$, which are of Fredholm for all $0 < \mu \leq \delta$.

(2)
Suppose the indices  are non negative. Then:
\begin{gather*}
0 \;\longrightarrow \; {\frak L}^2_{k+1}(X)  
\;\stackrel{d_{\mu}}{\longrightarrow} \; 
{\frak L}^2_k(X; \Lambda^1)  
  \stackrel{d^+_{\mu}}{\longrightarrow} \; 
L^2_{k-1}(X; \Lambda^2_+) \;
\longrightarrow \; 0.
\end{gather*}
is a family   of  Fredholm complexes of the same indices  for all $0 \leq \mu \leq \delta$.
\end{prop}
{\em Proof:} 
(2) follows from (1) with corollary $1.2$.
Noitce that $H^0_{\mu} =0$ always hold over non compact manifolds.

For (1),  we review the construction of the weighted Sobolev spaces
for convenience. For the details of the analysis, we refer to [K1].

{\bf Step 1:}
Let $(M,g)$ be a  closed Riemannian $3$ manifold,
and denote   the product metric 
by $g + dt$ on $ M \times {\mathbb R}$.   
By use of  formally $L^2$ adjoint operator,
we obtain the  elliptic operator
$P= d^*  \oplus  d^+ $ from $ \Lambda^1(M \times {\mathbb R})$ to
$ \Lambda^0(M \times {\mathbb R}) \oplus \Lambda^2_+(M \times {\mathbb R})$.
One can canonically identify:
$$\Lambda^1(M \times {\mathbb R}) = p^*(\Lambda^1(M)) \oplus p^*(\Lambda^0(M)), \quad
\Lambda^2_+(M \times {\mathbb R}) = p^*(\Lambda^1(M))$$
where $p: M \times {\mathbb R} \mapsto M$ is the projection, and the isomorphisms are
given by:
$$u+vdt \leftrightarrow (u,v), \quad *_M u+ u \wedge dt \leftrightarrow u.$$
Then
$P: p^*(\Lambda^1(M) \oplus \Lambda^0(M))
\mapsto  p^*(\Lambda^1(M) \oplus \Lambda^0(M))$ is represented as:
$$P= - \frac{d}{dt} + \begin{pmatrix}
  *_M d &d \\
d^* & 0
\end{pmatrix}
\equiv - \frac{d}{dt} + Q$$
where  $Q$ is an elliptic self adjoint  differential operator
on $L^2(M ; \Lambda^1 \oplus \Lambda^0)$.

{\bf Step 2:}
 Let us fix a small and positive  $\delta >0 $. 
Then for 
  $0 < \mu \leq \delta$,  define:
$$\tau: M \times [0, \infty) \mapsto [0, \infty), \quad \tau(m,t) = \mu t.$$
Let $X$ be a cylindrical manifold whose 
end  is isometric to $M \times [0, \infty)$.
Then we fix the weight function $\tau :X \to [0, \infty)$ so that it coincides with $\tau(m,t) $ 
  on the end of $X$.
Then we  define the
 weighted Sobolev $k$ norms on $X$ by:
$$||u||_{(L^2_k)_{\mu}}= (\Sigma_{l \leq k} \int_X \exp(\tau)|\nabla^l u|^2 )^{\frac{1}{2}}.$$
We write by $(L^2_k)_{\mu}$ as 
 the space of the completion of $C^{\infty}_c(X)$ with respect to the norm,
  since the isomorphism class of the function spaces
is determined by $\mu>0$, rather than $\tau$ itself.

By this way, we obtain a filtration of the Sobolev spaces $\{(L^2_k)_{\mu}\}_{0 \leq \mu \leq \delta}$,
which satisfy the inclusions whenever $\mu' \leq \mu$:
$$C^{\infty}_c \ \subset  \ (L^2_k)_{\mu}  \ \subset  \ (L^2_k)_{\mu'} \ \subset  \ L^2_k$$
such that the inclusions $(L^2_k)_{\mu}  \subset  L^2_k$ are dense for all $0 \leq \mu \leq \delta$.

{\bf Step 3:}
Let us  introduce the  isometries:
$$I_{\mu}: L^2(X, \Lambda^*)
\mapsto (L^2)_{\mu}(X,  \Lambda^*)$$
by $I_{\mu}(u)= \exp(-\frac{\tau}{2})u$,
which induce the isomorphisms:
 $$I_{\mu} : L^2_k(X) \cong (L^2_k)_{\mu} (X).$$

Let $d^*_{\tau}$ be the $(L^2)_{\mu}$ adjoint operator so that 
$<u, d(v)>_{(L^2)_{\mu}} = <d^*_{\tau}(u), v>_{(L^2)_{\mu}}$ holds,
and put:
$$P_{\tau} = d^*_{\tau} \oplus d^+:
L^2_{k+1}(X: \Lambda^1) \to L^2_k(X: (\Lambda^0 \oplus \Lambda^2_+))$$
Then we have the following expression
on the end  $M \times [0, \infty)$:
$$
I_{\tau}^{-1} P_{\tau} I_{\tau}  = - \frac{d}{dt}  + 
 \begin{pmatrix}
  *_M d & d \\
d^* & -\frac{d \tau}{dt}
\end{pmatrix}
+ \frac{1}{2}\frac{d \tau}{dt} \equiv   
 - \frac{d}{dt}  + Q_{\mu}
 $$
The following lemma is well-known.
Theorem $1.1$ with  lemma $1.4$ below  finishes the proof of (1):
\begin{lem}
The AHS$_{\mu}$ complexes:
\begin{gather*}
0 \;\longrightarrow \; (L^2_{k+1})_{\mu}(X)  
\;\stackrel{d}{\longrightarrow} \; 
(L^2_k)_{\mu}(X; \Lambda^1)  
  \stackrel{d^+}{\longrightarrow} \; 
(L^2_{k-1})_{\mu}(X; \Lambda^2_+) \;
\longrightarrow \; 0.
\end{gather*}
are of  Fredholm with $H^0_{\mu}=0$  for all $0 < \mu \leq \delta$.
\end{lem}
{\em Proof:}
 [K1] computed  spectral behavior of $Q_{\mu}$ near $0$.
For convenience let us outline how to verify this.
It follows  from the straightforward calculations that:
$$Q_{\mu} = \begin{pmatrix}
*_M d  + \frac{\mu}{2}  & d \\
d^* & -\frac{\mu}{2}
\end{pmatrix} : L^2_{k+1}(M  , \Lambda^1  \oplus  \Lambda^0 )
\mapsto L^2_k(M  , \Lambda^1 \oplus  \Lambda^0  )$$
give  isomorphisms for all small $ \mu >0$.
In particular $I_{\tau}^{-1} P_{\tau} I_{\tau}$ and hence $P_{\tau}$
both give the invertible operators on the end.
This completes the proof.

\vspace{3mm}

\section{Analysis over Casson handles}
{\bf 2.A Casson handles:}
Casson handles are open smooth four manifolds with the attaching regions.
They are inductively constructed by taking  end connected sums 
and   obtained as their direct  limits.
Each  building block  is called a
 kinky handle which    is diffeomorphic to a finite number of 
  the end connected sums $\natural (S^1 \times D^3)$
 with  two attaching regions, where
 one is a tubular neighborhood
 of bund sums of Whitehead links (this is connected with the previous block), 
 and the other  is a disjoint union of the 
standard open subsets $S^1 \times D^2$ in 
 $\sharp S^1 \times S^2 = \partial (\natural S^1 \times D^3)$ 
 (this is connected with the next block). The number of end-connected sums 
 is exactly the  one of self-intersections of the immersed two handle.
 
 By construction, each Casson handle corresponds to the  infinite rooted trees with sign
on each edge.
We attach a Casson handle to the zero handle along
 the first stage attaching circle and denote it by $S = D^4 \cup CH$.
 This is a smooth  open four manifold.
 We refer to [K1] for detailed description on  Casson handles.
\vspace{3mm} \\
{\bf 2.A.2 Kinky handles as Riemannian manifolds:}
There are two  simplest Casson handles  with respect to  signs.
Both have $S^1 \times D^3$ as their building blocks, and 
 are given by taking end connected sums of infinitely many kinky handles
periodically.

We have introduced the complete Riemannian metrics on kinky handles in 
[K1] (section $2$).
Let us briefly describe their properties.

Let $W_1 \cong S^1 \times D^3$ be the simplest kinky handle. The Riemannian metric on $W_1$
has the properties:

(1) $W_1$ contains two disjoint Riemannian subspaces:
$$N \times (- \infty ,0] \sqcup M \times [0, \infty)$$
where $N$ and $M$   are cylindrical  three manifolds with their ends $\Sigma$, 
which are mutually isometric.

(2) $W_1$ contains another Riemannian subspace:
$$V  \times [0, \infty)$$
where $V$  is a non compact Riemannian three manifold with two ends
which are isometric to $\Sigma \times (- \infty, 0] \cup \Sigma \times [0, \infty)$.

\vspace{2mm}

Let us denote:
$$\tilde{W}_1 = W_1 \  \backslash \ \{ \ N \times (- \infty, -2) \sqcup M \times (2, \infty) \ \}.$$
By taking 
the end connected sum of $N \times\{-2\}$ with $M \times \{2\}$,
we obtain the cylindrical four manifold:
$$Y_1 \equiv \tilde{W}_1 \  / \ \{ \  N \times \{-2\} \sim M \times \{2\}  \ \}.$$
Let us describe the periodic cover of $Y$.
Let us prepare infinitely many copies of $\tilde{W}_1$ and assign indices as $\{ \tilde{W}_1^i\}_{ i \in {\mathbb N}}$.
Then we take the  end connected sums of 
$M^i \equiv M \times \{2\}$ in $\tilde{W}_1^i$ with $N^{i+1} \equiv N \times \{-2\}$ in $\tilde{W}_1^{i+1}$.
By this way, we obtain the half periodic Riemannian manifold:
$$\tilde{Y}^0_1 \equiv CH({\mathbb N}) \equiv \tilde{W}^0_1 \cup \dots\dots \cup \tilde{W}_1^i \cup_{M^i \cong N^{i+1}}  \tilde{W}_1^{i+1} \cup \dots$$
which is the simplest periodic Casson handle, where
the attaching region lies in the boundary of $\tilde{W}_1^0$.

Notice that  
 the periodic Riemannian manifold:
$$\tilde{Y}_1 \equiv CH({\mathbb Z}) \equiv  \dots \cup \tilde{W}_1^i \cup_{M^i \cong N^{i+1}}  \tilde{W}_1^{i+1} \cup \dots$$
is obtained by use of  indices as $\{ \tilde{W}_1^i\}_{ i \in {\mathbb Z}}$,
which is ${\mathbb Z}$ covering of $Y_1$.
This space is less interesting from the view point of smooth structure,
since it is in fact diffeomorphic to the standard four disc.
However
its Riemannian structure 
 plays an important role in the analysis of Fourier-Laplace transform over
the period cover.
In fact we verify the following:
\begin{prop}
The AHS complex:
\begin{gather*}
0 \;\longrightarrow \; {\frak L}^2_{k+1}(\tilde{Y}_1; \Lambda^0)  
\;\stackrel{d}{\longrightarrow} \; 
{\frak L}^2_k(\tilde{Y}_1; \Lambda^1)  
  \stackrel{d^+}{\longrightarrow} \; 
L^2_{k-1}(\tilde{Y}_1; \Lambda^2_+) \;
\longrightarrow \; 0
\end{gather*}
is acyclic Fredholm.
\end{prop}
We verify this at the end of $2.B$ below.
\vspace{3mm} \\
{\bf 2.B Fourier-Laplace transform:}
Let $X$ be a complete Riemannian manifold, and $D : C^{\infty}_c(E) \to C^{\infty}_c(F)$ be a differential operator 
between vector bundles over $ X$.

Let $\tilde{X} $ be a periodic cover of $X$ with the group ${\mathbb Z}$.
Then both   $E,F$ and $D$ lift canonically 
as the invariant  operator:
$$\tilde{D}:  C^{\infty}_c (\tilde{E}) \to C^{\infty}_c (\tilde{F})$$
 where $\tilde{E} \to \tilde{X}$ is the natural lift, 
 equipped with the shift isomorphism
 $T: \tilde{E} \cong \tilde{E}$ which corresponds to $1 \in {\mathbb Z}$.
 
 Let $\tilde{W} \subset \tilde{X}$ be a fundamental domain with respect to ${\mathbb Z}$ action,
with the boundary $\partial \tilde{W}= N \cup M$.

\begin{defn}
Let us take any $\psi \in  C^{\infty}_c(\tilde{E}) $ and $z \in {\mathbb C}^*$.
The Fourier Laplace transform of $\psi$ is given  by:
$$\hat{\psi}_z(\quad) = \Sigma_{n=-\infty}^{\infty} z^n(T ^n \psi)(\quad)$$
over the restriction $\hat{\psi}| \tilde{W}$, which determines 
 a section  of the vector bundle:
$$E'  \equiv [ \tilde{E} \otimes_{{\mathbb R}} {\mathbb C} ]/ {\mathbb Z} \mapsto X \times {\mathbb C}^*$$
where $1 \in {\mathbb  Z}$ sends $(\rho, \lambda) \in \tilde{E} \otimes_{{\mathbb R}} {\mathbb C}$
to $(T \rho, z\lambda)$.
\end{defn}
{\em Remark 2.1:}
$E'(z)$ is a family of bundles over $X$ with $E'(1)=E$.
Every $E'(z)$ is isomorphic to $E$ ([K1] section $4$).

The {\em Fourier Laplace inversion formula} is given as follows;
for any  smooth section $\hat{\eta} \in C^{\infty}_c(E')$ with
 $\hat{\eta}_z \in  C_c^{\infty}(E'(z))$,  let us take the lift and restrict 
it  on $\tilde{W}$.
Then for any $s \in (0, \infty)$ and  $x \in \tilde{W}$,
$$T^n \eta(x) \equiv \frac{1}{2 \pi i} \int_{|z| =s} z^{-n} \hat{\eta}_z (\pi(x)) \frac{dz}{z}$$
defines a smooth section over $\tilde{E} \to \tilde{X}$, where
 $\pi: \tilde{W} \to X$ is the projection.
These are converses each other.

Let us define  the family of differential operators over  $E'(z)$
 by:
$$\hat{D}_z \hat{\psi}_z \equiv (\hat{D \psi})_z$$
passing  through the Fourier Laplace transform.

Suppose $X=Y$ is a cylindrical manifold and $\tau$ be a weight function on $Y$ with weight $\mu>0$.
Then the weight function canonically  extends on the periodic cover $\tilde{Y}$, and hence 
Sobolev spaces  $(L^2_k)_{\mu}(\tilde{Y})$  are obtained.
Let us extend $\tilde{D}$ and $D_z$ over $(L^2_k)_{\mu}(\tilde{E})$
and $(L^2_k)_{\mu}(E'(z))$ respectively.

\begin{lem}[K1]
$\tilde{D}$ is invertible over  $\tilde{Y}$, if
$D_z$ are invertible for all $z \in C(1) = \{z \in {\mathbb C} : |z|=1\}$.
\end{lem}
The assumption is satisfied, which we will explain in $2.B.2$ below.
\vspace{3mm}  \\
{\bf 2.B.2 Excision analysis:}
Let us introduce  some new  analytic method to bridge various functional spaces
over different spaces, which makes it convenient to perform excision process.

Let us prepare two complete Riemannian manifolds $X$ and $Y$
which satisfy the following conditions:

(1) There are open sub manifolds $A,X_0 \subset X$ and $B,Y_0 \subset Y$
such that $X'_0 \equiv X \backslash A$ and $Y_0' \equiv Y \backslash B$ 
are both manifolds with the product ends near  boundary:
$$X_0 \backslash X'_0 \cong [-3,-2] \times N, \quad  Y_0 \backslash Y_0' \cong [2,3] \times M$$
which are isometric $M \cong N$ mutually.

(2) The Dirac operators $D_X$ and $D_Y$ are equipped over $X$ and $Y$,
such that the isomorphism holds:
$$D_X| [-3,-2] \times N = D_Y|[2,3] \times M.$$

Let us denote the end connected sum:
$$Z = X_0 \cup Y_0 \  / \ \{ \   [-3,-2] \times N \sim [2,3] \times M \  \}$$
and denote by $D_Z$ as the induced Dirac operator over $Z$.

Let us consider $L^2_k(Y_0)_0$ and  regard it as a closed subspace
in $L^2_k(Z)$.

\begin{defn}
The orthogonal complement of 
$L^2_k(Y_0)_0$  in  $L^2_k(X)$ is given by:
$$L^2_k(Y_0)_0^{\perp} = \{ w \in L^2_k(X): <w, w'>=0 \text{ for any } w' \in L^2_k(Y_0)_0 \cap L^2_k(X)\}.$$
\end{defn}

\begin{lem}
$L^2_k(Y_0)_0^{\perp} $ splits into two closed  linear subspaces:
$$L^2_k(Y_0)_0^{\perp} = H_1 \oplus H_2$$
with
supp $H_1 \subset X_0$ and supp $H_2 \subset A$,
 such that
the orthogonal decomposition:
$$L^2_k(Z) = L^2_k(Y_0)_0 \oplus  H_1 $$
holds.
\end{lem}
{\em Proof:}
Any $u \in L^2_k(Z)$ can be expressed as a union 
$u_1+u_2 \in L^2_k(X_0)_0 + L^2_k(Y_0)_0$.
Let us take $v \in L^2_k(Y_0)_0^{\perp}$, and  verify that it defines a
continuous linear functional:
$$F_v: L^2_k(Z) \to {\mathbb R}, \qquad 
F_v(u) = <v, u_1> _{L^2_k(X)}.$$
Let us check its well-definedness.
Choose two different decompositions  $u =u_1+u_2 =u_1'+u_2'$.
Then $u_1-u_1' = u_2'-u_2 \in  L(X_0)_0 \cap  L^2_k(Y_0)_0$, and so:
$$<v, u_1-u_1' > _{L^2_k(X_0)_0} = <v, u_2'-u_2> _{L^2_k(X_0)_0}=0.$$
Thus we may regard $F_v =v_1 \in L^2_k(X_0)_0^* = L^2_k(X_0)_0 \subset L^2_k(Z)$.
This assignment gives the closed linear subspace
$H_1 \subset L^2_k(X_0)_0$.

Notice that    $X \backslash \{ [2,3] \times N\}$ splits into 
 the disjoint union 
$  X_0' \coprod A$ respectively.
 Now any $u \in L^2_k(X)$ can be expressed as a union: 
$$u_1+u_2+ u_3 \in L^2_k(X_0)_0 + L^2_k(Y_0)_0 \cap L^2_k(X_0)_0 + L^2_k(A)_0.$$
$<v, u_2>=0$ holds for $v \in L^2_k(Y_0)_0^{\perp}$.
Let us define another continuous  linear functional:
$$G_v: L^2_k(X) \to {\mathbb R}, \qquad 
G_v(u) = <v -v_1, , u> _{L^2_k(X)} =  <v -v_1, , u_3> _{L^2_k(X)} $$
So $G_v =v_2 \in L^2_k(A_0)_0^* = L^2_k(A_0)_0 \subset L^2_k(X)$,
which gives another closed linear subspace $H_2 \subset L^2_k(A_0)_0$.

In total  we have the decomposition:
 $$v  =v_1+ v_2 \in  L^2_k(X_0)_0  \oplus L^2_k(A)_0.$$
which  gives the closed linear subspaces $H_1 \oplus H_2 \subset L^2_k(X)$
with the desired properties.
This completes the proof.

\begin{prop}
Suppose  both  $D_Y$ 
and $D_X$  are of Fredholm.
Then $D_Z$  has  closed range.
\end{prop}
{\em Proof:}
{\bf Step 1:}
Suppose contrary. Then the spectrum of $D_Z$ accumulates near $0$.
Let us choose 
an orthonormal  sequence $\{u_i\}_i \subset L^2_{k+1}(Z)$ with $||u_i||_{L^2_{k+1}} =1$
and $||D(u_i)||_{L^2_k} \to 0$, where the spectra of  $u_i$ lie within $(- \lambda_i , \lambda_i)$ with $\lambda_i \to 0$.
Let us decompose $u_i = u_i^1 + u_i^2 \in L^2_{k+1}(Y_0)_0  \oplus H_1$ in lemma $2.1$.

We claim that there is positive $\epsilon >0$ with the uniform lower bounds:
$$\inf_i  \{ ||u_i^1||_{L^2_{k+1}}, ||u_i^2||_{L^2_{k+1}} \} > \epsilon.$$
Suppose contrary, and assume $||u_i^1||_{L^2_{k+1}} \to 0$.

Then both $||u_i^2||_{L^2_{k+1}} \to 1$ and 
$||D(u_i^2)||_{L^2_k} \to 0$ hold, while 
 support of $u_i^2 $ lie inside $X_0 \subset Z$.
 Moreover the inner products satisfy asymptotic orthogonality:
 $$<u^2_i, u^2_j>_{L^2_{k+1}} \to0 $$ holds as $i,j \to \infty$.
 
 If we regard $u_i^2 \in L^2_{k+1}(X)$, then there is $ m \geq 0$ 
 and $a^2_1+ \dots +a^2_m =1$ such that the uniform lower bound should hold:
 $$||\Sigma_{i=1}^m a_i D(u_i^2) ||_{L^2_k} \geq \delta >0$$
 since $D$ is of Fredholm over $X$.
However it follows from the assumption that
$||\Sigma_{i=1}^m a_i D(u_i^2) ||^2_{L^2_k} \to0$ should hold
as $i,j \to \infty$, which is a contradiction.

We can argue  another case $||u_i^2||_{L^2_{k+1}} \to 0$ by the same way.

This verifies the claim.

{\bf Step 2:}
Each $u_i$ are smooth since their spectra are small.
Let us decompose:
$$D(u_i) = w_i^1 \oplus w_i^2 \in  L^2_{k+1}(Y_0)_0  \oplus H_1$$ in lemma $2.2$.
Both $w_i^1$ and  $w_i^2$ converge to $0$ in $L^2_{k+1}$ as $i \to \infty$.

We claim that  $D(u_i^l) $ converge  to $0$ in $L^2_{k+1}$ for $l=1,2$.
Notice the equality
 $D(u_i^1) + D(u_i^2) =w_i^1+w_i^2$, and so:
$$\delta_i \equiv D(u_i^1) - w_i^2 = w_i^1 - D(u_i^2) \in L^2_{k+1}(Z).$$

For any $v   \in  C^{\infty}_c(Y_0)$, 
\begin{align*}
|<\delta_i, v >_{L^2_{k+1} }| & =|  < w_i^1, v> -<D(u_i^2), v> | \\
& =| <w_i^1, v> -  <u_i^2, D^*(v)>|  \\
& = |<w_i^1,v> | \leq ||w_i^1||_{L^2_{k+1}} ||v||_{L_{k+1}^2 } \to 0.
\end{align*}
 Thus  if we regard $\delta_i \in (L^2_{k+1}(Y_0)_0)^* \cong L^2_{k+1}(Y_0)_0 $, then it
converges to zero. So 
 $D(u_i^1) = \delta_i + w_i^2$ 
converge to zero in $L^2_{k+1}(Z)$.
However  this leads to a contradiction by arguing as in    step $1$.
 This completes the proof.
 \vspace{3mm}

 Let us apply proposition $2.2$ to analysis of the 
 parametrized elliptic operators  over $Y$.
 Let us consider
 the periodic cover:
  $$\tilde{Y} = \dots \tilde{W}^0 \cup \tilde{W}^1 \cup \dots$$
 where $\tilde{W}^i$ are the copies of the same $\tilde{W}$.
Then we obtain the  parametrized differential operators for $z \in C(1) \subset {\mathbb C}$:
 $$D_z: L^2_{k+1}(E(z)) \to  L^2_k(F(z)).$$

\begin{cor} Suppose $D: L^2_{k+1}(E) \cong L^2_k(F)$
gives an isomorphism. Then
 $D_z$
 have closed range
 for all $z \in C(1)$.
 \end{cor}
{\em Proof:}
 It is true for $z=1$ by the assumption.
  Let us denote $Y= \tilde{W}/ N(-2) \sim M(2)$, and 
$H, H_1  \subset \L^2_{k+1}(Y; E(z))$ be closed subspaces as in lemma $2.2$
satisfying:

(1) $H \oplus H_1 = L^2_{k+1}(Y; E(z))$ and

(2) $H =L^2_{k+1}( \tilde{W}; E(z))_0$
and Supp$H_1 \subset N \times [-2, 0]  \cup M \times [0,2]$.
\vspace{2mm} \\
Then one can identify  $H \subset L^2_{k+1}(Y; E)$.

For $u \in H_1$,  one can associate $u' \in L^2_{k+1}(Y; E)$ by: 
$$u'|N \times [-2, 0] =u, \quad u'|M \times [0,2] = z^{-1}u.$$
Let us denote by $\tilde{H}_1 \subset  L^2_{k+1}(Y; E)$ the corresponding 
subspace to $H_1$.
Since this assignment is isometric,  $\tilde{H}_1$ is also a closed subspace
of $L^2_{k+1}(Y; E)$.
Since $ D_i(\tilde{H}_1)$ is closed by the assumption,
 it follows that
$D_i(H_1) $
  is also a closed subspace.
  Then we can follow the proof of proposition $2.2$  and obtain the result.
  This completes the proof.

\vspace{3mm}

Let $Y$ be the cylindrical manifold given by the end connected sum of the kinky handles in $2.A.2$,
and let 
$\Lambda^*(z)$ for $z \in  C(1) = \{ z: |z| =1 \}$ 
be the twisted differential forms  in $2.B$.
The differentials canonically extend over the twisted
AHS complexes.
Moreover the weighted $L^2$ inner products are also induced from the one at $z=1$.
In [K1], we have verified the following:
\begin{lem}[K1] 
There is small $\delta >0$ so that the twisted AHS$_{\mu}$ complexes:
\begin{gather*}
0 \;\longrightarrow \; L^2_{k+1}(Y; \Lambda^0(z))  
\;\stackrel{d}{\longrightarrow} \; 
L^2_k(Y; \Lambda^1(z))  
  \stackrel{d^+}{\longrightarrow} \; 
L^2_{k-1}(Y; \Lambda^2_+(z)) \;
\longrightarrow \; 0
\end{gather*}
are acyclic for all $0 < \mu \leq \delta$.
\end{lem}

{\em Proof of proposition $2.1$:}
Combining with lemma $2.1$ and $2.3$, it follows that AHS$_{\mu}$ complexes
are acyclic Fredholm over the periodic cover $\tilde{Y}$
for all $0 <\mu \leq \delta$.
Then proposition $2.1$ follows from proposition $1.1$
(see also proof of proposition $1.2$).
This completes the proof.

\vspace{3mm}

The above method is the basis for the analysis over the higher stage Casson handles.
\vspace{2mm}  \\
{\bf 2.C Riemannian manifolds of the second stage:}
Let $(W_1,N,M)$ be the  simplest Riemannian kinky handle as described in $2.A.2$.
Let us recall  $Y_1 = \tilde{W}_1/ \{ N \sim M\}$ and its half periodic cover:
$$\tilde{Y}_1= \tilde{W}_1^0 \cup_{M^0 \sim N^1} 
\dots  \cup \tilde{W}_1^{i-1} \cup_{M^{i-1} \sim N^i}  \tilde{W}_1^i  \cup \dots $$
equipped with the attaching region $N^0 \subset  \tilde{W}_1^0 \subset \tilde{Y}_1$.
The half periodic Casson handle can be expressed  by the half real line $ T_1^0 =
 {\mathbb R}_{\geq 0}\equiv {\mathbb R}_0$
assigned with the same sign on each edge.
Let $T_2^0$ be the  tree: 
 $$T_2^0 = {\mathbb R}_0  \cup_{n \in {\mathbb N}} {\mathbb R}_0$$
 which is obtained from $T_1^0$ attached with the infinite number of the same  $(T^0_1)'$ at the root with 
each integers ${\mathbb N} \subset {\mathbb R}_0$.

Let us describe the corresponding Riemannian-Casson handle $CH(T_2) =Y_2$.
Let $(W_2, N,M_1,M_2)$ be a kinky handle with $2$ kinks. 
$W_2$ is obtained by the end connected sum of two copies of $W_1$ along the disks on the boundary:
$$W_2 = W_1  \ \natural  \ W_1.$$
The Riemannian structure on $W_2$  satisfies the following properties:

(1) $W_2$ contains three disjoint Riemannian subspaces:
$$N \times (- \infty,0] \cup M_1 \times [0, \infty) \cup M_2 \times [0,\infty)$$
where $N$ and $M_1, M_2$ are cylindrical manifolds with  their  ends $\Sigma$, which are mutually isometric.

(2) $W_2$ contains another Riemannian subspace:
$$V_2  \times [0, \infty)$$
where $V_2$  is a non compact Riemannian three manifold with three ends
which are  isometric to three disjoint union of $\Sigma \times [0, - \infty)$.

\vspace{2mm}

Let us denote:
$$\tilde{W}_2= W_2 \backslash \{ N \times (- \infty, -2) \cup M_1 \times (2, \infty) \cup M_2 \times (2,\infty) \}$$
and 
 attach $CH({\mathbb N})$ in $2.A.2$ by 
identifying $M_2 \times \{2\}$ with $N^0 \times \{-2\}$ in $\tilde{W}^0_1 \subset CH({\mathbb N})$:
$$CW_2 = \tilde{W}_2 \cup_{M_2 \sim N^0} CH({\mathbb N}).$$
If we take the end connected sum along  $N \times\{-2\}$ with $M_1 \times \{2\}$ as:
$$Y_2  \equiv  CW_2 \  / \   \{ N \times \{-2\} \sim M_1 \times \{2\} \}$$
then $Y_2$ is a complete Riemannian manifold without boundary.

Let us consider 
 its half  periodic cover:
$$CH(T_2^0) = \tilde{Y}^0_2 =  CW_2^0 \cup_{M_1^0 \sim N^1}  \dots \cup CW_2^i \cup_{M_1^i \cong N^{i+1}}  CW_2^{i+1} \cup \dots$$
where $\{ CW_2^i\}_{ i \in {\mathbb N}}$ are the infinite number of the copies of $CW_2$.
This is the Casson handle corresponding to $T_2^0$ with the attaching region $N^0 \subset CW_2^0$.

 As above, we denote the periodic cover of $Y_2$ by:
$$\tilde{Y}_2 =   \dots \cup CW_2^i \cup_{M_1^i \cong N^{i+1}}  CW_2^{i+1} \cup \dots$$

Combining proposition $2.1$ and $2.2$,
AHS$_{\mu}$ complex  has closed range for $0 < \mu \leq \delta$ over $Y_2$.
It has been verified to be acyclic by use of an asymptotic method in [K1].
\vspace{3mm} \\
{\bf 2.D Casson handles of higher stages:}
By use of  Fourier-Laplace transform with
 the  parallel  argument to use Fourier-Laplace transform
with acyclicity of AHS$_{\mu}$ complexes over $Y_2$ above, 
it follows that the  AHS$_{\mu}$  complexes over $\tilde{Y}_2$ are also acyclic Fredholm.

  Let us consider the third stages. Let:
  $$CH(T_2^0) = \tilde{Y}_2^0= CW_2^0 \cup  \dots \cup CW_2^i \cup_{M_1^i \cong N^{i+1}}  CW_2^{i+1} \cup \dots$$
  be the second stage of the Casson handle, 
  and put:
$$CW_3 = \tilde{W}_2 \cup_{M_2 \sim N^0} CH(T^0_2), \quad 
Y_3  \equiv  CW_3 \  / \   \{ N \times \{-2\} \sim M_1 \times \{2\} \}.$$
Let us consider:
$T_2^0 = {\mathbb R}_+ \cup_{n  \in {\mathbb N}} {\mathbb R}_+$, and put:
$$T_3 = {\mathbb R} \cup_{n \in {\mathbb Z}} T_2^0.$$
The induced periodic cover is the Casson handle which correspond to $T_3$:
$$\tilde{Y}_3 =CH(T_3).$$
The same argument as above verifies that 
 the  AHS$_{\mu}$  complexes over  both $Y_3$ and $\tilde{Y}_3$ are  acyclic Fredholm.

By this way, we obtain the complete Riemannian Casson handles $Y_N$ so that 
 the  AHS$_{\mu}$  complexes over  both $Y_N$ and $\tilde{Y}_N$ are  acyclic Fredholm.

So far we have described the rooted trees with at most trivalent branch.
The above construction works for the rooted trees with more branches.
In [K1], we have introduced a class of homogeneous  trees of bounded type.
In our notation, the half periodic trees $T^0_k $ are expressed as  $T^0_{k+1} = T^0_{2,2, \dots , 2,1}$, where $2$ appear   $k$ times.
Let $n_1, \dots, n_k \in \{ 1,2, \dots \}$ be positive integers. 
Then  using kinky handles with $n_j$ kinks, one has  a natural extension, and gets 
the homogeneous tree of bounded type
 $T^0_{n_1, \dots, n_k,1} $ which is a rooted  infinite tree and admits 
 the corresponding Riemannian-Casson handle  $CH(T^0_{n_1, \dots,n_k, 1})$.

By iterating the previous process, one can verify that
the  AHS$_{\mu}$  complexes over:
   $$CH(T_{(n_1, \dots, n_k,1)} )\equiv \tilde{Y}_{(n_1, \dots, n_k,1)}$$ 
    are all acyclic Fredholm.

For practical application,  one considers open four manifolds given by
the $0$ handle attatched with Casson handles.
Recall that $ l(S^2 \times S^2) \backslash $ pt is homotopy equivalent to
the wedges of $S^2$, and  also that 
it  is described  pictorically by $l$ disjoint unions of Hopf links with $0$ framings,
where each $S^1$ component corresponds to the attaching region of the Casson handle.
So there is a diffeomorphsm:
$$ l(S^2 \times S^2) \backslash \text{ pt } 
 \cong D^4 \natural \cup_{l=1}^{2l} (D^2 \times D^2).$$

Let $T_1, \dots, T_{2l}$ be signed homogeneous trees of bounded type.
Let us equip with a complete Riemannian metric on $D^4$ with boundary $M \times \{0\}$,
which 
contains $M \times [0,1]$ isometrically.
By the end connected sum, one obtains the Riemannian-Casson handle:
 $$S \equiv  D^4 \natural \cup_{l=1}^{2l} CH(T_l).$$

 \begin{thm}[K1]  Let $S = D^4 \cup^{2l}_{j=1} CH(T_j)$ 
 be the Riemannian-Casson handle 
whose 
trees are  homogeneous of bounded type.  Then
there is a complete Riemannian metric $g$ 
of bounded geometry on $S$ and  positive $\delta >0$ 
so that  the  bounded complexes:
$$
0 \;\rightarrow \; (L^2_{k+1})_{\mu}(S,g ) 
\;\stackrel{d}{\longrightarrow} \; 
(L^2_k)_{\mu}((S,g); \Lambda^1 ) 
  \stackrel{d^+}{\longrightarrow} \; 
(L^2_{k-1})_{\mu}((S,g); \Lambda^2_+ ) \;\rightarrow \; 0$$
are   Fredholm for all the weight $0 < \mu \leq \delta$
with their cohomology groups:
\begin{align*}
 & H^0_{\mu} =0, \quad H^1_{\mu}=0, \\
& l=\dim H^2_+(S:{\mathbb R})  \leq  \dim  H^2_{\mu} \leq 2l = \dim H^2(S:{\mathbb R}).
\end{align*}
\end{thm}

Now combining theorem $2.1$  with corollary $1.2$, we obtain:
\begin{cor} 
The  AHS complex  over $S$:
\begin{gather*}
0 \;\longrightarrow \; {\frak L}^2_{k+1}(S)  
\;\stackrel{d}{\longrightarrow} \; 
{\frak L}^2_k(S; \Lambda^1)  
  \stackrel{d^+}{\longrightarrow} \; 
L^2_{k-1}(S; \Lambda^2_+) \;
\longrightarrow \; 0.
\end{gather*}
is  of Fredholm with the same index above.
\end{cor}

\section{Trivializing  at infinity}
{\bf 3.A Tree like structure on Riemannian manifolds:}
Let $X$ be a complete Riemannian  four manifold of bounded geometry,
and $E \to X$ be an $SO(3)$ bundle over $X$.
Each $x \in X$ admits local chart onto $\delta_0 >0$ ball in ${\mathbb R}^4$.

A  covering $\{U_i\}_i$ on a
Riemannian manifold $X$ is  {\em bounded}, 
if  (1) $\sup_i $ diam $U_i < \infty$ and (2) $\sup_i  \sharp \{ j: U_i \cap U_j \ne \phi\} < \infty$ hold.
Later on we assume that $U_i$ are diffeomorphic to the disks, which 
actually  gives no extra conditions  on the existence of bounded coverings,
when  $X$ is of bounded geometry.

Let $T^*$ be a connected and rooted tree with the root $*$, and 
 introduce the canonical tree metric on it.
For any vertex $v \in T^*$, let $g(v) \subset T^*$ be the set of vertices
which lie on the geodesic from the root to $v$.

\begin{defn}
(1)
A Riemannian manifold $X$ admits  $T^* \times {\mathbb N} $ covering,
if there exists a bounded covering $\{U_i\}_i$  
and a one to one correspondence $I: T^* \times {\mathbb N}  \cong {\mathbb N}$
such that  $I$ satisfies: 
\begin{align*}
& U_{I(k,i)} \cap U_{I(l,j)}= \phi \text{ for } d(l,k) \geq 2, \\
&
\sup_{k,m,i,l} 
 \{  |i-m| : U_{I(l,i)} \cap U_{I(k,m)}   \ne \phi \} < \infty  .
\end{align*}

(2)
$T^* \times {\mathbb N}$ covering over  $X$ is
trivial at infinity, if  for any $C$, 
there are  $C_0, k_0, j_0$ such that 
 for any $v \in T^*$ with  $|v|  \geq k_0$ and $j \geq j_0$,
  $C$ neighborhood of: 
 $$ \cup_{ v' \in g(v)}  U_{I(v,j)} \cup_{ 0 \leq i \leq j}  U_{I(v,i)} $$
 is  $\pi_1$-null in its  $C_0$ neighborhood. 
\end{defn}

Let us see particular cases:

(1)
 A Riemannian manifold $X$ admits  ${\mathbb N}$ covering,
if there exists a  bounded 
covering $\{U_i\}_i$ and  a one to one correspondence $I: {\mathbb N} \cong {\mathbb N}$
with:
$$   \sup_{i,j} \{|i-j| : U_{I(i)} \cap U_{I(j)} \ne \phi \} < \infty.$$
{\em Example 3.1:}
Any cylindrical manifolds admit  ${\mathbb N}$ covering.
\vspace{3mm}

(2)
A Riemannian manifold $X$ admits  ${\mathbb N}^2$ covering,
if there exists a bounded 
covering $\{U_i\}_i$ and a one to one correspondence $I: {\mathbb N}^2 \cong {\mathbb N}$
with:
$$ U_{I(k,i)} \cap U_{I(l,j)}= \phi \text{ for } |l-k| \geq 2, \quad
\sup_{k,m,l,i} 
 \{ |i-m| : U_{I(l,i)} \cap U_{I(k,m)}   \ne \phi \} < \infty .$$

\begin{lem}
Let  $S \equiv  D^4 \natural \cup_{l=1}^m CH(T_l)$
be the Riemannian-Casson handle in $2.D$.
Then $S$ admit $T^* \times {\mathbb N}$ covering
which is trivial at infinity.

In particular 
 the half periodic Casson handles admit ${\mathbb N}^2$ covering
 which is trivial at infinity.
\end{lem}
{\em Proof:}
{\bf Step 1:}
Let us describe how the half periodic Casson handles admit ${\mathbb N}^2$ covering.
Let:
 $$S =  D^4  \ \natural  \ CH({\mathbb N}) =
 D^4  \ \natural  \ \tilde{W}_1^0 \ \natural \ \tilde{W}_1^1 \ \natural \ \tilde{W}_1^2 \natural \dots$$
where  $\tilde{W}_1^i$ contain two disjoint product ends 
$N \times [-1,0] \coprod M \times [0, 1]$, and 
 $D^4$  has one product end  $N \times [-1,0]$.
 On the other hand 
 $D^4$ and $\tilde{W}^i_1$ are manifolds with boundary which contain 
 cylindrical ends $B \times [0, \infty)$ and $V^i \times [0, \infty)$ respectively,
  where    the Riemannian three manifolds $V^i$ are all  mutually isometric
  with two boundary components.
  
 Let us decompose both $B $ and $V$ by  open discs
 $B = \cup_{j=1}^{m_0} U_j$ and  $V = \cup_{j=1}^{m_0} U'_j$.
Let $I: {\mathbb N}^2 \to S$ be a map which satisfies the following properties:

(1) The image 
$I(0, {\mathbb N }) \subset D^4$ 
is  a net such that $I(0, i \ m_0) \in B \times \{ i + i_0 \}$ for some $i_0$.

(2) the image
$I(k, {\mathbb N }) \subset \tilde{W}_1^k$ 
is a net such that $I(k, i \ m_0) \in  V^k \times \{ i + i_0\} $.
\\
It is immediate to see that $I$ satisfies the required conditions.

Let us consider the general case, and let $m$ be the number of the trees in the Casson handles.
For $m =1$, $S$ clearly admits $T_1^* \times {\mathbb N}$ covering
in a parallel way.

For $m \geq 2$, let $T^* $ be the rooted tree which attach all $T_l$ at their roots.
Then the  construction  of $T^* \times {\mathbb N}$ covering  can be reduced to the case $m=1$.

{\bf Step 2:}
Let us consider triviality at infinity.
For simplicity of the notation, we verify the half periodic case $S= D^4 \natural CH({\mathbb N})$ only.
Let:
$$CH({\mathbb N}) \equiv \tilde{W}^0_1 
\cup \dots\dots \cup \tilde{W}_1^k \cup_{M^k \cong N^{k+1}}  \tilde{W}_1^{k+1} \cup \dots$$
be the half periodic Riemannian-Casson handle equipped with 
${\mathbb N}^2$ covering.
The end of $\tilde{W}_1^k$ is isometric to $V^k \times [0, \infty)$.

It is enough to see that
for any $i, k$,
$$(\cup_{0 \leq l \leq k} V^l  ) \times \{i\}  \cup M_k $$
 are $\pi_1$-null in:
  $$(\cup_{0 \leq l \leq k+1} V^l ) \times \{i\} \cup (\tilde{W}^{k+1}_1 \backslash V^{k+1} \times [i+1, \infty)).$$
This follows from the construction of the kinky handles, where 
 the embedding of the solid torus 
$T^1 \equiv S^1 \times D^2 \hookrightarrow S^1 \times D^2 \equiv T^2$ is given 
by the Whitehead double, and hence $T^1$ is contractible in $T^2$.
This completes the proof.

\vspace{3mm}

Let $Y$ and $Z$ be two metric spaces. 
They are mutually {\em quasi-isometric}, if there is a map
$f: Y \to Z$ (not necessarily continuous) 
and  some $C,D$ so that (1)
$C$ neighborhood of  the image of $f$ covers $Z$, and 
(2) the inequalities hold:
$$C^{-1}d_Y(m,m') -D \leq d_Z(f(m),f(m')) \leq Cd_Y(m,m') +D$$
 for any $m,m'  \in Y$.

\begin{cor} 
Any  map $f: T^* \times{\mathbb N} \to S$ which satisfies $f(v,m) \in U_{I(v,m)}$
gives a   quasi-isometry:
$$T^* \times {\mathbb N} \cong S= D^4 \natural \cup_{l=1}^m CH(T_l)$$
with respect to 
$T^* \times {\mathbb N}$ structure 
in lemma $3.1$.
\end{cor}
{\em Proof:}
This follows from the construction of the Riemannian metrics and $T^* \times {\mathbb N}$
structure equipped above.
This completes the proof.
\vspace{3mm} \\
{\bf 3.A.2: Example:}
Let us start from the following:
\begin{lem}
Hyperbolic space ${\bf H}^4$ does not admit
 quasi-isometry with $T^* \times {\mathbb N}$.
\end{lem}
{\em Proof:}
The asymptotic dimension is a numerical invariant 
for metric spaces, which is preserved under quasi-isometry  ([Gr]).
So  their dimensions should coincide with each other,
 if they could admit quasi-isometry.
However 
$\text{as-dim } {\bf H}^4$
is in fact equal to $4$, while
 $\text{ as-dim }  (T^* \times {\mathbb N}) $
 is $2$.
This completes the proof.

\vspace{3mm}

M.Tsukamoto observed the following:
\begin{lem}
Let ${\bf H}^4$ be as above. Then for any small $r \in [0, 1)$, there exists 
an $L^2$ ASD connection $A_r$ with 
$$||F_{A_r}||L^2({\bf H}^4) =r.$$

In particular it should not necessarily integer.
\end{lem}
{\em Proof:}
Let $A$ be a non trivial ASD connection 
with $ e = ||F_A||L^2(D) >0$
over the unit disc $D \subset {\mathbb R}^4$.
Then for any $0< a \leq e$, there exists some $0 < r \leq 1$ with
$ a = ||F_A||L^2(D_r)$, where $D_r \subset D$ is $r$-disc.

Notice that the ASD condition and $L^2$ norm are both preserved
under conformal change of the metrics.
Let us equip with Poincar\'e metric on $D_r$ which is isometric to ${\bf H}^4$.
Then the restriction $A|D_r$ with the metric is the desired one.
\vspace{3mm} \\
{\bf 3.B Trivializing near infinity:}
Let us recall 
the  topology of the end of $S$:
\begin{thm}[Fr]
The end of
   Casson handle  $S = D^4 \cup^{2l}_{j=1} CH(T_j)$ 
   admits a topological color $\cong S^3 \times [0, \infty)$.
\end{thm}
In particular $S$ is simply connected and simply connected at infinity.

We verify the following:
\begin{prop}
Let $X$ be a complete Riemannian manifold
of bounded geometry
 equipped with 
 $T^* \times {\mathbb N}$ covering.
Let $E \to X$ be a bundle on $X$ 
and $A$ be  an   ASD connection over $E$
with $||F_A||L^2(X) < \infty$.

 If the covering is trivial  at infinity, then 
 $A$ is approximated by a  compactly supported  smooth connection,
after gauge transformation.

\end{prop}
{\em Proof:} Notice that if
$A$ is an ASD connection whose curvature $F_A$ is in $L^2$, then
 local $L^2$ norms are sufficiently small 
$||F_A||_{B_{\delta_0}(x)} < \epsilon$ near infinity $x \in X$.

 In order to verify this, one uses the following.
 \begin{sublem}[U1]
Let  $A$ be an ASD connection over $E$, and choose a local  trivialization
$E| B_{2\delta_0}(x) \cong {\mathbb R}^3$.
Then there exist $\epsilon_0 >0$ and $C_k$ for  $k \geq 3$
such that
if  $||F_A||_{L^2(B_{2\delta_0}(x))} < \epsilon_0 $ holds,
then there exists a gauge transform
$g \in $ Aut $E|B_{2\delta_0}(x)$ with:
$$g^*(A) = d+ A', \quad ||A'||_{L^2_k(B_{\delta_0}(x))} \leq  C_k ||F_A||_{L^2(B_{2\delta_0}(x))}.$$
\end{sublem}
This is Uhlenbeck's theorem based on the construction of the Coulomb gauges.
In particular $A$ is smooth by the Sobolev embedding.
 
 \vspace{3mm}

{\em Proof of proposition:}
We split the proof  into $5$ steps, where we will not use
triviality at infinity on $X$ until step $4$.

{\bf Step 1:}
Let  $\{U_i\}_i$ be a bounded  covering,
where each $U_i$ is diffeomorphic to the disc.
Let us choose any frame at $m \in U_i$, and construct  the local  frame: 
$$ \psi_i : E|U_i \cong U_i \times {\mathbb R}^3.$$
 by use  of the  parallel transport with respect to $A$.
 
 We claim that 
for any $\epsilon >0$,  there exists $i_0$ such that for any $i \geq i_0$,
the restriction can be expressed as:
 $$A|U_i = d +a_i, \quad ||a_i||L^2_k(U_i) < \epsilon$$ 
  with respect to the  trivialization over $U_i$. 
In fact  
$||F_A|| L^2(U_i) < \epsilon$ hold for all large $i>>1$ by the assumption. Then
 it follows from sublemma $3.1$ that there is a local trivialization
  $ \varphi_i: E|U_i \cong U_i \times {\mathbb R}^3$ such that 
the above estimates hold with respect to $\varphi_i$.
One may assume that $\psi_i$ and $\varphi_i$ coincide at $m \in U_i$.
Recall that the parallel transport is given uniquely by the ODE:
$$\nabla_{\dot{x}(t)} \xi = \Sigma \{ \frac{d \xi^{\lambda}(x(t))}{dt} 
+ \Sigma \Gamma_{\mu, i}^{\lambda}(x(t)) \frac{dx^i}{dt} \xi^{\mu}(x(t))\} \ e_{\lambda}$$
where $\nabla e_{\lambda} = \Sigma \ \omega_{\lambda}^{\mu} e_{\mu}$, 
$A= (\omega_{\lambda}^{\mu})$ with
$\omega_{\lambda}^{\mu} = \Sigma \ \Gamma_{\lambda, i}^{\mu} dx^i$.
In particular  we may assume that $\Gamma_{\lambda,i}^{\mu}$
have sufficiently small norms in $L^2_k$  with respect to $\varphi_i$.
This verifies the claim, since the parallel transport is independent of choice of local coordinates.

{\bf Step 2:}
Recall that if $x_1(t)$ and $x_2(t)$ are two homotopic paths with the same end points, 
then the difference of their  parallel transports $\xi^1$ and $\xi^2$  with the same initial vector $\xi$
can be estimated by $l^1$ norms of the curvature of $A$ over 
any surface which span $x_1 \cup x_2$.

As a test case, let us consider a simple situation with
$X= U_1 \cup U_2 $, and trivializations
 $ \varphi_i: E|U_i \cong U_i \times {\mathbb R}^3$
   are given over $U_i$.
   Let us  choose any $m \in U_1 \cap U_2$.
If  $g_{12} = \varphi_1^{-1} \varphi_2(m)  \in SO(3)$ is sufficiently near the identity, 
then small modification of $\varphi_2|U_1 \cap U_2$ gives the  extension of the trivialization $\varphi_1$ 
over $U_1 \cup U_2$. Notice that this is the case under the situation in step $1$, 
by change of  the trivialization over $U_2$ by constant if necessarily.

Next let us consider $X= U_1 \cup U_2 \cup U_3$ with the  trivializations
 $ \varphi_i: E|U_i \cong U_i \times {\mathbb R}^3$
   over $U_i$ as above.
   If $X= U_1 \cap U_2 \cap U_3 \ne \phi$, then
   we can replace their trivializations over $U_i$ by constants 
   so that small modification of their new  trivializations over $U_i$  gives the global one over $X$.
   Suppose $X= U_1 \cap U_2 \cap U_3 = \phi$, and 
   extensions of the  trivializations over $U_1 \cup U_2$ and $U_2 \cup U_3$  are given as above.
   Let us  choose any $m \in U_1 \cap U_3$ and consider    $g_{13} = \varphi_1^{-1} \varphi_3(m)  \in SO(3)$.
 If  $ \varphi_1^{-1} \varphi_3$ takes value sufficiently near the constant $g_{13}$, 
  one can modify $\varphi_3$ on $U_1\cap U_3$ so that 
  the result of the local trivializations
  give flat structure over $X$. Notice that if $g_{13}$ is away from constant, then 
  the flat structure may not be able to extend over $U_1 \cup U_2 \cup U_3 $.

  {\bf Step 3:}
  Let us consider a case 
 with many number of the coverings
  $X= U_1 \cup \dots \cup U_i \cup \dots  $, and suppose:
  
  (1)  
 there is some $s_0$   such that  $|i-j| \leq s_0$ hold whenever
$U_i \cap U_j \ne \phi$, 

(2) 
 the curvature of $A$ has small $L^2$ norms over each $U_i$ as in step $1$.
\\
Then one can construct flat structure over $X$  inductively
by use of the method in step $2$ as below.

Let us draw a bi-Lipschitz  line $l : [1, \infty) \to X$ 
so that $l(i) \in U_i$ hold.  Choose any frame at $l(1) \in U_1$ and 
fix the trivializations on all $l(i) \in U_i$ by use of the parallel transport.
By use of the trivializations at $l(i)$, one obtains the local trivializations over $U_i$ as in step $1$.
Let us extend 
 flat structure inductively.
Suppose it is given over   $U^{m-1} \equiv  U_1 \cup \dots \cup U_{m-1}$, and let us extend it over $U^m$.
If $U_i \cap U_m \ne \phi$, then $i \geq m-s_0$ must hold.

As a simple case, let us assume moreover:
\\
 (3) there is some constant $C$ such that $\cup_{i \leq j \leq i+s_0} U_j $ are contractible   
 in $C$ neighborhood of them for all $i$.

In this case we obtain the global trivialization over $X$, 
since 
the transition functions $\varphi_i^{-1} \varphi_m$ over $U_i \cap U_m$ must take
values sufficiently near the identity,
as we noticed at the first paragraph of step $2$, and
 the trivialization can be extended as in step $2$.

Now let us remove the condition (3), and construct flat structure over $X$.
Let us fix a constant $C >> \sup_i \text{ diam }  \{ U_i \cup \dots \cup U_{i+s_0}  \}$. 
Suppose $x \in U_i \cap U_m \ne \phi$ for some $ i \leq m$.
Let us consider  the loop $l_{i,m}$ 
given by the union of $l|[i, m]$ with another line between $l(i)$ and $l(m)$ in $U_i \cup U_m$.

If $l_{i,m}$ is contractible in $C$ neighborhood of $l(m)$, then 
$\varphi_i^{-1} \varphi_m$ take values near the identity, and 
one can modify $\varphi|U_i \cap U_m$ slightly so that 
trivialization is extended over $U_i \cap U_m$ to $U_i \cup U_m$.

If $l_{i,m}$ is not contractible in $C$ neighborhood of $l(m)$, then 
one can modify $\varphi|U_i \cap U_m$ slightly so that 
 flat structure  is extended over $U_i \cap U_m$ to $U_i \cup U_m$ as in step $2$.

Next choose another $i' \leq m$ with $U_{i'} \cap U_m \ne \phi$.
If $U_i \cap U_{i'} \cap U_m = \phi$, then we can extend the flat structure over $U_{i'} \cup U_m$ as above.
Suppose $U_i \cap U_{i'} \cap U_m \ne \phi$.
If $l_{i',m}$ is homotopic to $l_{i,m}$  in $C$ neighborhood of $l(m)$, then 
$\varphi_{i'}^{-1} \varphi_i$ take values near the identity, and so 
one can modify $\varphi|U_{i'} \cap U_m$ slightly,
 preserving it over  $U_i \cap U_{i'} \cap U_m $, 
so that one can extend the 
trivialization over $U_{i'} \cap U_m$ to $ U_m$.
If $l_{i',m}$ is not homotopic to $l_{i,m}$  in $C$ neighborhood of $l(m)$, 
then 
one can modify $\varphi|U_{i'} \cap U_m$ slightly so that one can 
extend flat structure  over $U_{i'} \cup U_m$ to $U_m$ as in step $2$.

One can repeat this process at most $s_0$ times, and obtain the extension of the flat structure over $U^m$.
This finishes the induction step.

Notice that 
 if  $U_i \cap U_j \ne \phi$ could happen with large $|i-j|$, 
then the corresponding transition functions will vary very large.
This is the key aspect where we have introduced uniformity of the multiplicity of coverings.

It is a basic fact that for a compact $X$, two bundles over $X$ are mutually homotopic, if
their transition functions are sufficiently near. In fact there is some $n$ such that 
there are bundle surjections $\varphi : X \times {\mathbb R}^n \to E$ 
and $\varphi' : X \times {\mathbb  R}^n \to E'$
by use of the trivializations over $U_i$ respectively.
Since their transition functions are mutually near, the homotopy:
$$x \in X \to t \text{ Ker } \varphi (x) +(1-t) \text{ Ker } \varphi'(x)$$
gives a family of bundles over $X$. In particular Ker $\varphi$ and Ker $\varphi'$ are
mutually isomorphic. Then $E \cong (X \times {\mathbb R}^n) / \text{ Ker } \varphi$
is isomorphic to $E'$ (see [A] p $29$).

Suppose $X$ admits  ${\mathbb N}$ covering.
Let us check that  $A$ can be approximated by a 
smooth connection which is flat at infinity.

Let us  choose $j_0$ with $I(j) \geq i_0$ for all $j \geq j_0$, where $i_0 $ is in step $1$. 
Let us take a trivialization over $U_{I(j_0)}$ as in step $1$.

By step $3$, one obtains the flat structure on the end of $E$ with the estimates:
$$A|\cup_{j_0 \leq j } U_{I(j)} = d+a, \quad ||a||L^2_k(\cup_{j \geq j_0}  U_{I(j)}) < \epsilon.$$

{\bf Step 4:}
Suppose $X$ admits 
${\mathbb N} \times {\mathbb N} $ covering
which is trivial near infinity.
Let us choose $j_0$ and $k_0$ 
such that   $L^2$ norm of $F_A$ is less than small $\epsilon >0$ over:
$$\cup_{ j \geq j_0, k \geq 0}  U_{I(j,k)}  \cup_{j \geq 0, k \geq k_0} U_{I(j,k)} .$$

Let us put the set:
 $$L(j_0,k_0)  = \{ (j,k_0) :  0 \leq j \leq j_0\} \cup \{(j_0,k):  0 \leq k \leq k_0\}$$
 and $L(j_0,k_0,s_0) = \cup_{j_0-s_0 \leq j \leq j_0+s_0} \cup_{k_0-s_0 \leq k \leq k_0+s_0} L(j,k)$. 
 Correspondingly, we put:
 $$U(j_0,k_0,s_0) = \cup_{j_0-s_0 \leq j \leq j_0+s_0} \cup_{k_0 -s_0\leq k \leq k_0+s_0} U_{I(j,k)}.$$ 
 Let us fix $s_0' >>s_0>>0$ such that $U(j_0,k_0,s_0)$ is $\pi_1$-null in $U(j_0,k_0,s_0')$.

 Let $\bar{L}(j_0,k_0)$ be the connected line in $[0, \infty) \times [0, \infty)$ which connects 
 neighbor points by each straight line in $L(j_0,k_0)$.
Let us choose a bi-Lipschitz  line:
$$l: \bar{L}(j_0, k_0) \to X, \quad l(j,k) \in U_{I(j,k)}$$ so that it connects the points 
 in $U_{I(j_0,k)} $ and $U_{I(j_0,k+1)}$ for all $ 0 \leq  k \leq k_0-1$,
and the points  in $U_{I(j,k_0)}  $ and $U_{I(j+1,k_0)}$ for all   $  0 \leq j \leq  j_0-1 $.

Let us choose a frame at a point in   $U_{I(j_0,0)}$ as in step $1$, and 
extend the trivialization along the line by parallel transport.

By step $1$, extend the trivialization  over  $U_{I(j,k)} $  for all $(i,j) \in L(j_0,k_0,s_0')$.

By step $3$,    extend the flat structure inductively over $U(j_0,k_0,s_0') $ along the line.
By restriction, we obtain the trivialization over $U(j_0,k_0,s_0) $:
$$A|U(j_0,k_0,s_0)= d+a, \quad ||a||L^2_k(U(j_0,k_0,s_0)) < \epsilon.$$
 It follows from sub lemma $3.1$ that $A$ can be approximated by the trivial connection
 over $\cup_{(j,k) \in L(j_0,k_0,s_0)} U_{I(j,k)}$ by use of cut off function on the region.

{\bf Step 5:}
Let us consider  the case of 
$ T^* \times  {\mathbb N} $ covering which is trivial at infinity.
The idea is similar to step $4$.
Let us put $T^*_k = \{  v \in T^* : d(v,*) \leq  k\}$ with $\partial T^*_k = \{  v \in T^* : d(v,*) = k\}$.
Then we choose $j_0$ and $k_0$ 
such that 
  $L^2$ norm of $F_A$ is less than $\epsilon$ over:
 $$T(j_0,k_0)  = \{ (j,v) :  0 \leq j \leq j_0, v \in \partial T^*_{k_0} \} \cup \{(j_0,v ):  v \in T^*_{k_0} \}$$
 Let us put  $T(j_0,k_0,s_0) = \cup_{j_0-s_0 \leq j \leq j_0+s_0} \cup_{k_0-s_0 \leq k \leq k_0+s_0} T(j,k)$,
 and:
 $$U(j_0,k_0,s_0) = \cup_{j_0 -s_0\leq j \leq j_0+s_0} \cup_{k_0 -s_0\leq k \leq k_0+s_0} U_{I(j,k)}.$$

Let $\bar{T}(j_0,k_0) $ be the connected tree in $T^* \times [0, \infty)$
which connects 
 neighbor points by each straight line in $T(j_0,k_0)$.
Let us choose   a bi-Lipschitz map
$l: \bar{T}(j_0,k_0) \to X$
with $l(j,k) \in U_{I(j,k)}$ for all $(j,k) \in T(j_0,k_0)$.
Then as in step $4$, we choose trivialization by parallel transport
along  $\bar{T}(j_0,k_0)$.
 The rest process is the same as step $4$.

This completes the proof.

\vspace{3mm}

\begin{cor} Suppose $X$ admits $T^* \times {\mathbb N}$ structure which is trivial at infinity.

Then
  $p_1(A) = \frac{1}{4 \pi^2} \int_Y tr(F_A \wedge F_A)$ is an integer.
\end{cor}
{\em Proof:}
Let us choose the trivialization as above.
Then there is a  compact subset $K \subset X$, 
so that the bundle over $X \backslash K$ is trivial and
 $A|X \backslash K = d+ a$ with 
$  ||a||L^2_k(X \backslash K) < \epsilon$ is sufficiently small.

By cut off, $A$ is approximated by a compactly supported smooth connection.
One may assume that the boundary of $K$ is a smooth submanifold of codimension  $1$.
Let us consider the double of $K$, and 
extend the approximated connection over it, where it is trivial over the extra $K$.
Then  $p_1$ must take integer value over the double,
which is equal to the $p_1$ value  for the approximated connection over $X$.
In particular $p_1(A)$ itself also takes integer value.
 
This completes the proof.

\vspace{3mm}

\section{Local analysis of the moduli space of ASD connections}
Let $X$ be a non compact smooth four manifold which is 
 simply connected and simply connected at infinity.
Let $g$ be a complete Riemannian metric of bounded geometry 
so that 
  injectivity radius is uniformly bounded from below by 
 a positive constant $\epsilon >0$,
  and   the curvature operator satisfies uniform bound from above  as
 sup$_{x \in X} |\nabla^l R| < \infty$  for any $l \geq 0$.

Let $E \to X$ be  an $SO(3)$ vector bundle which 
is trivial over $X \backslash K$ for 
 some compact subset $K$, 
  equipped with a fixed connection  $\nabla_0$ which is trivial near infinity
 with respect to the trivialization. 
  $E$ is determined by $w_2(E) \in H^2(X:{\mathbb Z}_2)$
 and $p_1(E) \in {\mathbb Z}$.
Let $P$ be  the corresponding principal $G$ bundle
with $P \times_G {\mathbb R}^3 =E$, and 
 put the adjoint bundle by:
 $$ Ad(P) = P \times_G {\frak G}$$ where 
$\frak G$ is the Lie algebra of $G$.
\vspace{3mm} \\
{\bf 4.A Connection spaces:}
Let $A$ be an ASD connection over $E$
such that $a_0 \equiv A - \nabla_0  \in L^2_k(X ; Ad(P) \otimes \Lambda^1)$
for any $k \geq 0$.
We call such a connection as an $L^2$ ASD connection.

The Atiyah-Hitchin-Singer complex (AHS complex) is given by:
\begin{gather*}
0 \; @>>>\;  
C_c^{\infty}(X; Ad(P))  \; @> d_A >>\; 
    C_c^{\infty}(X; Ad(P) \otimes  \Lambda^1 )
\; @> d^+_A >>\;
C_c^{\infty}(X; Ad(P) \otimes \Lambda^2_+) \; @>>>\; 0
\end{gather*}
where $d^+_A= (1+*) \circ d_A$.

Let us introduce the corresponding functional spaces, which
heavily depend on choice of $A$:
\begin{defn}
The functional spaces 
 $ {\frak L}_l(A)$
are given 
by the maximal extension of the domain
$C^{\infty}_c(X; Ad(P) \otimes \Lambda^*)$ 
 with \ their norms:
$$||u||^2_{{\frak L}_{k+1}^2(A)} =||d_A(u)||^2_{L^2_k}, \quad
||w||^2_{{\frak L}_k^2(A)} = ||d_A(u)||^2_{L^2_k} + ||d^+_Aw'||^2_{L^2_{k-1}}
$$
with respect to the orthogonal decomposition $w=d_A(u)+w'$ as:
$$L^2_k(X; Ad(P) \otimes \Lambda^1) 
= d_A({\frak L}_{k+1}^2(X; Ad(P))) \oplus d_A({\frak L}_{k+1}^2(X; Ad(P) ))^{\perp}.$$ 
\end{defn}
{\em Remark 4.1:}
In general 
one cannot  integrate $d_A$ in order to obtain
gauge group actions in a straightforward way, since 
 $g^{-1} ag $ will not  be  in 
${\frak L}_k(A)$
 for any gauge group $g$ and $a \in {\frak L}_k(A)$.
 If we try to approximate  $a $ by compactly supported smooth forms,
their  $C^0$ norms may grow  unboundedly,
even though $d^+_A(a)$  keep bounded $L^2_{k-1}$ norms.
Such situation will happen when the spectrum of the Laplacian on $1$ form 
contain the continuous part near zero.
\vspace{3mm} \\
{\bf 4.B ASD connections:}
Let $A $ be an $L^2$ ASD connection.
Then the isomorphisms  of the Hilbert spaces hold  by lemma $1.2$:
 \begin{align*}
& d_A : {\frak L}_{k+1} (A) \cong  d_A({\frak L}_{k+1}(A)) 
\subset L^2_k(X: Ad(P) \otimes  \Lambda^1) , \\
& {\frak L}_k(A) \cong d_A({\frak L}_{k+1}(A)) \oplus d^+_A({\frak L}_k(A) ) \\
& \qquad \quad
  \subset L^2_k(X: Ad(P) \otimes \Lambda^1) \oplus L^2_{k-1}(X: Ad(P) \otimes \Lambda^2_+).
\end{align*}

\begin{lem}
 $(\text{ Ker } d^+_A)^{\perp} \subset L^2_k(X: Ad(P) \otimes \Lambda^1)$
is dense in $(\text{ Ker } d^+_A)^{\perp} \subset {\frak L}_k(A)$.
  \end{lem}
  {\em Proof:}
   $(\text{ Ker } d^+_A)^{\perp} \subset {\frak L}_k(A)$ is isomorphic to
  $d^+_A({\frak L}_k(A))$, and the closure of 
 $d^+_A((\text{ Ker } d^+_A)^{\perp} ) \subset L^2_{k-1} (X: Ad(P) \otimes \Lambda^2_+)$
 is $d^+_A({\frak L}_k(A))$.
 
 This completes the proof.

\vspace{3mm}

Let $A , A'$ be   two $L^2$ ASD connections
with:
 $$||a||L^2_k(X) \equiv ||A-A'||L^2_k(X)  < \infty$$ 
 
 \begin{lem}
 Suppose $[A] =[A'] \in  {\frak L}_k(A)$.
 Then 
 
 (1) $A_t =A + t(A-A')$ is a parametrized ASD connections.
 
 (2) The curvatures are all the same $F_{A_t} = F_A$.
 \end{lem}
 {\em Proof:}
 The assumption implies $d^+_A(a) = (a \wedge a)^+ =0$ hold.
 In particular for any $t \in {\mathbb R}$, 
 $A_t =A+ ta$ gives a family of ASD connections, since the equalities hold:
 $$F^+_{A_t} = F^+_A + t d^+_A(a) + t^2(a \wedge a)^+=0.$$
 
 On the other hand $||F_{A_t}||L^2(X)$ must be constant,
 which implies $d_A(a)=a \wedge a=0$ must hold.
 So 
 $F_A=F_{A_t}$ hold.
 This completes the proof.
 \vspace{3mm} \\
  {\bf 4.C Reguralization:}
    Let $A$ be an $L^2$ ASD connection and   consider: 
$$F^+:  A+ L^2_k(X:  Ad(P) \otimes  \Lambda^1) \to  L^2_{k-1}(X: Ad(P) \otimes \Lambda^2_+).$$
In general it would be impossible to extend it as a continuous functional 
from ${\frak L}_k(A)$, since kernel of $d^+_A$ will affect to determine output.
Even when we restrict it on the orthogonal complement of the kernel, 
still spectra near zero will affect when  the image of $d^+_A$ is not closed.

Let us fix an $L^2$ ASD connection $A$ and $\epsilon >0$, and
consider the spectral decomposition:
$$d^+_A \circ (d^+_A)^* = \Delta_A = \int_0^{\infty} \lambda^2 dE(\lambda)$$
on $L^2(X: Ad(P) \otimes \Lambda^2_+)$.
Let $f: [0,\infty) \to {\mathbb R}$ be a smooth function,
and put:
\begin{align*}
& P_f= \int_0^{\infty} f(\lambda) d E(\lambda), \\
&
Q_f(a) = (d^+_A)^* \Delta_A^{-1} P_f( d^+_A(a)).
\end{align*}
Then we  introduce a deformation of smooth functionals:
\begin{align*}
& F^+_f :  L^2_k(X: Ad(P) \otimes \Lambda^1)  \to L^2_{k-1}(X: Ad(P) \otimes \Lambda^2_+) \\
& d_A^+(a) +  (Q_f(a) \wedge Q_f(a))^+
\equiv d^+_A(a) + B_f(a,a)
\end{align*}

As a particular case, let
 $f_{\epsilon}: [0,\infty) \to [0,1]$ be a smooth function
with: 
$$f_{\epsilon}(\lambda) =
\begin{cases} 
 0  & \lambda \in [0,\epsilon] \\
 1  & \lambda \geq 2\epsilon
\end{cases}$$

 \begin{lem}
 (1) Suppose $d^+_A: L^2_k \to L^2_{k-1}$ has closed range.
 If we restrict $F^+_{f_{\epsilon}}$ over $(\text{ ker } d^+_A)^{\perp}$, then
it coincides with the standard self-dual curvature functional
 for all sufficiently small $\epsilon >0$.

 (2) 
  $ Q_{f_{\epsilon}}: {\frak L}_k(A) \to L^2_k $ gives a bounded linear functional:
  $$ ||Q_{f_{\epsilon}} (\quad ) ||L^2_k \leq 
C_k \epsilon^{-1}   ||\quad ||{\frak L}_k(A)  $$
  where $C_k$ depends only on $k$.
  
 (3)
$$F_{f_{\epsilon}}^+: {\frak L}_k(A) \to L^2_{k-1}(X: Ad(P) \otimes \Lambda^2_+)$$
defines a uniformly Lipschitz functional.
In particular the differential at $A$ is isometric onto its image
over $(\text{ Ker } d^+_A)^{\perp}$.
\end{lem}
So  the operator norm $||(d^+_A)^{-1}|| =1$ holds,
when $d^+_A$ is an isomorphism on $(\text{ Ker } d^+_A)^{\perp}$.

\vspace{3mm}

{\em Proof of lemma $4.3$:}
(1) if we choose $\epsilon < \lambda_1$ smaller than the first eigenvalue 
of $\Delta_A$ on $ L^2(X: Ad(P) \otimes \Lambda^2_+) $, 
then $ Q_{f_{\epsilon}}(a)=a$ hold for all $a \in (\text{ ker } d^+_A)^{\perp}$.

For (2), notice the estimate:
\begin{align*}
<(d^+_A)^*  \Delta_A^{-1} (P_{f_{\epsilon}}(b)), (d^+_A)^*  \Delta_A^{-1} (P_{f_{\epsilon}}(b))>_{L^2} 
& =  <P_{f_{\epsilon}}(b),  \Delta_A^{-1} (P_{f_{\epsilon}}(b)) >_{L^2} \\
& \leq \epsilon^{-1} <b,b>_{L^2}.
\end{align*}
Then it 
 follows from the estimates below:
\begin{align*}
||(d^+_A)^*  \Delta_A^{-1} P_{f_{\epsilon}}( d^+_A(a))||_{L^2_k}  
   & \leq C_k \epsilon^{-1} 
 ||P_{f_{\epsilon}}( d^+_A(a))||_{L^2_{k-1}} \\
&  \leq C_k \epsilon^{-1} ||d^+_A(a)||_{L^2_{k-1}} = C_k \epsilon^{-1}||a||_{{\frak L}_k(A)}.
 \end{align*}

(3) Let $a_t = a+ tb \in {\frak L}_k(A)$, 
and consider the difference:
$$
F_{f_{\epsilon}}^+(a+tb)- F_{f_{\epsilon}}^+(a)  =
t d^+_A(b) + B_{f_{\epsilon}}(a,a) + t(B(a,b)+B_{f_{\epsilon}}(b,a) )+ t^2 B_{f_{\epsilon}}(b,b).$$
The following estimates hold by the H\"older  estimate and  (2):
$$||B_{f_{\epsilon}}(a,b)||_{L^2_{k-1} } \leq C_k  ||Q_{f_{\epsilon}}(a)||_{L^2_k} ||Q_{f_{\epsilon}}(b)||_{L^2_k}
\leq C_k \epsilon^{-1} ||a||_{{\frak L}_k(A)} ||b||_{{\frak L}_k(A)}$$
where the constant $C_k$ is independent of $a,b \in {\frak L}_k(A)$.

In particular the
 Lipschitz constant of the differential   $d^+_A$ at $A$ is $1$, by definition of 
 the norm on ${\frak L}_1(A)$.
This completes the proof.

\vspace{3mm}

\begin{defn} Let us fix $\epsilon >0$.
The regularization of  $F^+$ at $A$  is given by:
\begin{align*}
& \tilde{F}^+ : {\frak L}_k(A)  \to L^2_{k-1}(X: Ad(P) \otimes \Lambda^2_+), \\
& \tilde{F}^+ (A+a) = 
d_A^+(a) + (Q_{f_{\epsilon}}(a) \wedge Q_{f_{\epsilon}}(a))^+
\end{align*}
\end{defn}

The following holds by the implicit function theorem:
\begin{lem} 
(1) 
 If $0$ is  regular value of both:
  $$F_{f_{\epsilon}}^+, \ F^+ : L^2_k(X: Ad(P) \otimes \Lambda^1) 
\to   L^2_{k-1}(X: Ad(P) \otimes \Lambda^2_+)$$ 
 then there is a cobordism 
between regular smooth manifolds:
$$(F^+_{f_{\epsilon}})^{-1}(0), \ (F^+)^{-1}(0) \subset A + L^2_k(X: Ad(P) \otimes \Lambda^1) .$$

(2)  If $0$ is a regular value of $\tilde{F}^+$ over ${\frak L}_k(A)$, 
 then
$(\tilde{F}^+)^{-1} (0)\subset  {\frak L}_k(A) $
 is a regular smooth manifold.
 \end{lem}
{\em Proof:}
We only have to verify (1).
Let
 $f_{\epsilon}^t :  [0,\infty) \to [0,1]$ be a smooth path of functions
such that $f_{\epsilon}^1=f_{\epsilon}$ and $f_{\epsilon}^0 \equiv 1$.
Then one obtains a family of the smooth functional:
$$F^+_{f_{\epsilon}^t} :
  L^2_k(X: Ad(P) \otimes \Lambda^1)  \to L^2_{k-1}(X: Ad(P) \otimes \Lambda^2_+)$$
by $d_A^+(a) +  (Q_{f_{\epsilon}^t}(a) \wedge Q_{f_{\epsilon}^t}(a))^+ $.
Then  the conclusion follows by the inverse function theorem.
This completes the proof.

\vspace{3mm}

\begin{cor}
Let  $A$ be an $L^2$ ASD connection over  $X$, and assume
it  is regular so that $d^+_A$ is surjective.
Then  for any  $\epsilon>0$, there is $\delta >0$ 
which is independent of $X$ such that for 
any small perturbation 
$s$ with $s(b) \in L^2_{k-1}(X: Ad(P) \otimes \Lambda^2_+)$ with $||s(b)||_{L^2_{k-1}} \leq \delta$, 
there is a solution to the  perturbed  equation  in ${\frak L}_k(A)$:
$$(\tilde{F}^+ +s)(a)=0. $$
\end{cor}
{\em Proof:}
$a=0$ is the solution for $b =0$.
Then this follows from lemma $4.3(3)$ and the inverse function theorem.
This completes the proof.
\vspace{3mm} \\
{\em Remark 4.1:}
This uniformity plays a key role for our proof of theorem $0.2$.
The above property holds even if we replace 
$\tilde{F}^+$  by the linear functional $d^+_A$.
\vspace{3mm} \\
{\bf 4.D Slices and index formula:}
Let $U_{\epsilon} \subset {\frak L}_k(A)$ be $\epsilon$ neighborhood of $A$,
and put:
$$U_{\epsilon}^{\perp} \equiv \{ A+v \in {\frak L}_k(A) : v \in d_A({\frak L}_{k+1}(A))^{\perp} \cap U_{\epsilon} \}.$$

Let $A$ be an $L^2$ ASD connection, and consider 
$ (\text{ Ker } d^+_A)^{\perp} \subset L^2_k$.
Since closure of  im $d_A $ is contained in   Ker $d^+_A$,
there is an  embedding:
$$ (\text{ Ker } d^+_A)^{\perp} \subset {\frak L}_k(A)$$

Let  $V_{\epsilon} \subset L^2_k(X: Ad(P) \otimes \Lambda^1)$ be $\epsilon$ neighborhood, 
and put:
$$V_{\epsilon}^{\perp} = V_{\epsilon} \cap (\text{ Ker } d^+_A)^{\perp}.$$
Then the embedding
$V_{\epsilon}^{\perp} \hookrightarrow U_{\epsilon}^{\perp}$
has dense image by lemma $4.1$.

\begin{prop}  
Let 
 $S = D^4 \natural \cup_{i=1}^{2l} CH(T^i)$ is  the Rimannian-Casson handle  
 homogeneously of bounded type, 
and choose an $L^2$ ASD connection  $A$. 

Then
 for any $L^2$ ASD connection $A$,
 \begin{gather*}
0 \;\longrightarrow \; {\frak L}^2_{k+1}(A)
\;\stackrel{d_A}{\longrightarrow} \; 
{\frak L}^2_k(A)  
  \stackrel{d^+_A}{\longrightarrow} \; 
L^2_{k-1}(X;  Ad(P) \otimes \Lambda^2_+   ) \;
\longrightarrow \; 0
\end{gather*}
is of  Fredholm whose  index is larger than or equal to:
 $$2p_1(P) + 3   \text{ dim}H^2_{\mu} \ \geq \  2p_1 (P) + 3l$$
if the number is non negative.
\end{prop}
{\em Proof:}
This follows from corollary $1.2$ with the next lemma.
The index computation below uses the excision principle, or relative index theorem ([GL]).

\begin{lem} [K1] Let 
 $S $ and $A$ be as above. Then for each small $\mu >0$,
the AHS$_{\mu}$ complex:
\begin{equation}
\begin{align*}
0\; @>>>\; (L^2_{k+1})_{\mu}(S; Ad(P))  & \;   @> d_A >>\; 
(L^2_k)_{\mu}(S; Ad(P) \otimes \Lambda^1) \\
  & \;   @> d_A^+ >> 
(L^2_{k-1})_{\mu}(S; Ad(P) \otimes \Lambda^2_+)\; @>>>\; 0
\end{align*}
\end{equation}
is a Fredholm complex whose  index satisfies the bounds:
$$ 2p_1(P) +3l \ \leq \ 2p_1(P)  +3  \dim H^2_{\mu}  \  \leq  \ 2p_1(P) +6l.$$
\end{lem}
See also  theorem $2.1$ on the AHS complex without coefficient.

\begin{cor}
If   $d^+_A: {\frak L}_k(A) \to L^2_{k-1}$ is surjective,
then the 
local moduli spaces:
$${\frak M}(A)_{loc} = \{ A' =A+a \in U_{\epsilon}^{\perp} \subset  {\frak L}_k(A)
: \tilde{F}^+(A')=0  \}$$
is a regular manifold whose dimension is smaller than or equal to:
$$-2p_1(P) - 3 \dim H^2_{\mu}.$$

In particular if $2p_1(P) + 3l  >0$ is positive, then 
the regular manifold should have negative dimension, and hence 
$d^+_A$ can not be surjective in the case. 
\end{cor}

\section{Local perturbation and transversality}
{\bf 5.A Convergence process:}
Let us fix $\epsilon >0$ and $\delta>0$ in corollary $4.1$. 

Let $M$ be a closed smooth  four manifold, and choose a family of 
Riemannian metrics  $g_i$ of bounded geometry  on $M$ such that they converge 
on each compact subset
to a complete 
Riemannain metric $h$ on an open subset  $S \subset M$.

{\bf Process 1:}
Let $E$ be an $SO(3)$ bundle, and take a family of 
regular ASD connections $A_i$ with respect to $(M,g_i)$
 so that  $d^+_{A_i}$ are surjective and
 hence give isomorphisms on $(\text{ Ker }d^+_{A_i})^{\perp}$.

By taking a  subsequence of $\{A_i\}_i$, they  converge to an $L^2$ ASD conenction $A$
over $S$ on each compact subset.

Let $b \in C_c^{\infty}(S: Ad(P) \otimes \Lambda^2_+)$
with $||b|| \leq \delta$
be a smooth self dual $2$ form on $S$ with sufficiently small support.
It follows from corollary $4.1$ that 
there are family of solutions $A_i+a_i $ to the equations:
$$\tilde{F}^+_{A_i}(a_i)=b$$
 such that $a_i$ have uniformly bounded 
 ${\frak L}_k(A_i)$ norms over $(M, g_i)$.
 
 {\bf Process 2:}
 Let us verify that a subsequence of $\{a_i\}_i$ converge weakly
 to $a \in {\frak L}_k(A)$ over $S$.
 Let $K \subset S$ be a compact subset.
 
 Let us consider  a family of bounded linear maps:
$$\varphi_i \in L^2_{k-1}(K)_0^* \qquad
 u \to < u, d^+_{A_i}(a_i)>$$
where $||\varphi_i|| < \infty$ are uniformly bounded from above
independently of $K$.
Then there is $\varphi \in L^2_{k-1}(K)_0^*$ which is a 
 weak limit of $\varphi_i$. 
Let: 
$$\Psi: {\frak L}_k(A) \to L^2_{k-1}(K)_0^* \qquad a \to u \to <d^+_A(a),u>$$
be the bounded linear functional.

We claim that 
$\varphi_i$  are approximated by the image of $\Psi$.
In fact there are $a_i' \in L^2_{k-1}(M,g_i)$ which approximate $a_i$ by lemma $4.1$.
Then by use of the cut off function which are equal to one on $K$,
one may regard them as elements in $L^2_{k-1}(S,h)$ and hence  in ${\frak L}_k(A)$,
which verifies the claim.

{\bf Process 3:}
Next we claim  that 
$\Psi$ has closed range, if: 
$$d^+_A: {\frak L}_k(A) \to L^2_{k-1}(S:Ad(P) \otimes \Lambda^2_+)$$ has finite codimenson.
This follows from the following abstract argument.
Let $H$ be a Hilbert space and $V$ be a finite dimensional vector space.
Let $L \subset H \oplus V$ be a closed linear subspace, 
and consider $d: H \to L^*$ by the same way as above.
Let us see  the image of $d$ is closed.
Then it follows from this abstract property that $\Psi$ has closed range.

Let $P: H \oplus V \to H$ be the projection.
Since $d(h)(l) =<h,l> =<h, P(l)>$ hold, it is enough to see that $P(L) \subset H$ is closed.
One may assume $L \cap V=0$, since  $P(L \cap V)=0$ holds.

So suppose $P(l_i)$ converge to $h \in H$.
Let us decompose $l_i =h_i \oplus v_i$.
Then $h_i $ converge to $h$.
Suppose $v_i$ are bounded sequence in $V$.
Then by finite dimensionality, a subsequence converge to $v \in V$, and hence
$h+ v \in L$. Assume $v_i$ could be unbounded.
Then by rescaling by constants so that $||l_i||=1$ with $h_i \to 0$
and so $||v_i|| \to 1$.  Then a subsequence converge to $v \in L$ which 
contradicts to the assumption, and we are done.

{\bf Process 4:}
Then $\varphi_i$ lie in the image of $\Psi$, and 
 there is some $a \in {\frak L}_k(A)$ with: 
$$\lim_i \varphi_i(u) = <u, d^+_A(a)>_{L^2_{k-1}(S)}$$ hold for all
 $u \in L^2_{k-1}(K)_0$.
 Notice that for a closed linear subspace $L \subset H$ in a Hilbert space,
 if a sequence $a_i \in L$ weakly converge to some $a \in H$, then $a\in L$ holds,
 since $<a,v >=\lim_i <a_i, v>=0$ hold for all $v \in L^{\perp}$.

  $||a||{\frak L}_k(A)$
  is uniformly bounded independent of choice of $K$, and
 denote it  by $a_K$.

  Let us choose an exhaustion of $S$ by compact subsets:
 $$K_0 \subset K_1\subset \dots \subset S$$
 and choose the corresponding $a_i \equiv a_{K_i} \in {\frak L}_k(A)$
 which consists of uniformly bounded sequence.
 Again there is a weak limit $a = w-\lim_i a_i \in {\frak L}_k(A)$.

We claim that $a$ solves the equation:
$$\tilde{F}^+(a)=b.$$
Let $k^i(x,y)$  and $k(x,y)$ be the smooth kernels of $\Delta_{A_i}^{-1} P_{f_{\epsilon}}(\Delta_{A_i})$
and $\Delta_A^{-1} P_{f_{\epsilon}}(\Delta_A)$
over $(M,g_i)$ and $(S,g)$ respectively.
Then $k^i$ converge to $k$ smoothly on each compact subset in $S \times S$.
It follows from the Sobolev estimate  that $\tilde{F}^+(a_i)$ weakly converge to 
$\tilde{F}^+(a)$ so that
the equality holds:
$$<u, \tilde{F}^+(a)> = <u, b>$$
for any $u \in L^2_{k-1}(S: Ad(P) \otimes \Lambda^2_+)$.
Since $\tilde{F}^+(a)$ itself is in $L^2_{k-1}$,
it is equal to $b $.
This verifies  the claim.
\vspace{3mm} \\
{\bf 5.B Transversality:}
In order to obtain regular moduli spaces, 
we have to perform  transversality argument.
Uhlenbeck's  metric perturbation method does not seem to work for our case,
since the functional spaces we have introduced, 
may change their structure heavily under small perturbation of metrics.

Another approach by 
 holonomy perturbation surely does not change the structure of our functional spaces,
since perturbation is local. However they  are given by the equivalent classes of connections,
which can not determine holonomy.

For our purpose of proof of theorem $0.2$,
we do not require gauge invariant perturbations.
This makes the situation quite simple, and we take a simplest way 
just to perturb the self dual $2$ forms directly.

Let $D \subset S$ be a small disk,
and   put the perturbation space $B$
$$B =L^2_{k-1}(D: Ad(P) \otimes \Lambda^2_+)_0$$
 with the inclusion $i: B \hookrightarrow L^2_{k-1}(S: Ad(P) \otimes \Lambda^2_+)$.

\begin{lem} The  functional:
$$\tilde{F}^+ +i: {\frak L}_k(A) \times B \to L^2_{k-1}(X: Ad(P) \otimes \Lambda^2_+)$$
gives the  surjective differential at $(A,0)$.
\end{lem}
{\em Proof:}
Elements of the cokenel of the image $d^+_A({\frak L}_k(A)) \subset L^2_{k-1}$ 
 can be assumed to satisfy the equation
 $(d^+_A)^*(u)=0$, since the image is the closure of $d^+_A(L^2_k(X: Ad(P) \otimes \Lambda^1))$.
 
 So elements in  the cokenel of $d (\tilde{F}^++i)|_{(A,0)}: {\frak L}_k(A) \times B \to L^2_{k-1}(X: Ad(P) \otimes \Lambda^2_+)$
 also satisfy  the equation $(d^+_A)^*(u)=0$.
 By unique continuation property,  the restriction $u|D$  does not vanish.
 On the other hand 
 one can choose   $b \in B $ with $<b, u> \ne 0$, which gives a contradiction.
 This completes the proof.
 
 \vspace{3mm}

Let $U_{\epsilon}^{\perp} \subset {\frak L}_k(A)$ be in $4.D$.
It follows from  lemma $5.1$  with the infinite dimensional inverse function theorem, 
that the map:
$$\tilde{F}^+ +i : U_{\epsilon}^{\perp}  \oplus B \to  L^2_{k-1}(S: Ad(P) \otimes \Lambda^2_+ )$$
has $0 \in U_{\epsilon}^{\perp}  \oplus B$  as a regular point. So the inverse:
$${\frak M} B \equiv (\tilde{F}^+ +i)^{-1}(0) \subset {\frak L}_k \oplus B$$
is the infinite dimensional Hilbert manifold near zero.

\begin{cor} There is an open neighborhood
$U \subset  U_{\epsilon}^{\perp} \oplus B$ and 
 a Baire set $\tilde{B} \subset B$ 
such that $\tilde{B}$ is the regular set over $U$:
$${\frak M}(A,b) = \{ a \in {\frak L}_k(A)   : (a,b) \in U , \ \tilde{F}^+(a) + b = 0 \}$$
is a  regular and finite dimensional smooth manifold
for any $b \in \tilde{B}$.
Its dimension is  equal to  the codimension of 
$d^+_A: {\frak L}_k(A) \to L^2_{k-1}(S: Ad(P) \otimes \Lambda^2_+)$,
which  is smaller than or equal to:
$$-2p_1(P) - 3 \dim H^2_{\mu}.$$
\end{cor}
{\em Proof:}
Let $\pi: {\frak M}B \to B$ be the projection, which is 
of Fredholm by proposition $4.1$.
Then the conclusion follows by 
 the Sard-Smale theorem and corollary $4.2$.
 
 This completes the proof.
 \vspace{3mm} \\
{\bf 5.C Proof of theorem $0.2$ and dimension counting:}
Let us recall the ideas of the argument in [K3].
Let $M$ be $K3$ surface, and choose an $SO(3)$ bundle $E$ over $M$
such that  the Donaldson's invariant does not vanish over $E$ ([Kr]).
So there are always ASD connections over $E$ and generic ASD moduli spaces
have $0$ dimensional.

Let us proceed by contradiction argument to see 
that Casson handles of bounded type cannot be embedded into $K3$ surface.
Suppose it could be, then by definition, 
Casson handles of homogeneously bounded type can also be embedded.
Let us denote   the Riemannian-Casson handle
by $(S,h)$.

Let us choose a family of generic metrics $g_i$ over $M$ such that 
they converge to $h$ on each compact subset on $S \subset M$.

Let us choose any ASD connections $A_i$ with respect to $(M,g_i)$, 
which converges to an $L^2$ ASD connection $A$ over $S$, and 
 fix  a trivialization near infinity so that $A = d+m$ with $||m||L^2_k < \infty$ holds.
So $A_i =d+m_i$ with $||m_i||(L_k^2)_{\text{loc}} < \infty$ on each compact subset of $S \subset M$
with respect to the trivialization.

Let us apply process $1$ in $5.A$.
Let $B$ be the Banach perturbation space which consists of local sections on $D \subset S$
in $5.B$,  and choose solutions 
$\tilde{F}^+_{A_i}(a_i)=b$  for $a_i \in {\frak L}_k(A_i)$ for generic $b \in \tilde{B}$ 
over ${\frak L}_k(A)$ in $5.B$.
 $d^+_A: {\frak L}_k(A) \to L^2_{k-1}(S: Ad(P) \otimes \Lambda^2_+)$ has finite codimension,
 and so $\{a_i\}_i$ converge weakly to $a \in {\frak L}_k(A)$ by processes $2-4$.
$a$  is a solution to the equation
$\tilde{F}^+_A(a)+b =0$, and so the space:
$${\frak M}(A,b) = \{ a \in {\frak L}_k(A)  \cap U : \tilde{F}^+_A(a) + b = 0 \}$$
would be a  non regular non empty smooth manifold, 
whose dimension is  smaller than or equal to
$-2p_1(P) - 3 \dim H^2_{\mu} \leq -2p_1(P) - 9$ by proposition  $4.1$.

Let us estimate its formal dimension.
Over  $K3$ surface or its logarithmic transforms, 
we consider  the case with $p_1=-6$ with $l = 3$.
 After this deformation process of metrics, one obtains another bundle whose 
  absolute value of the  first Pontryagin number strictly decreases
(see the argument also  in [K1]).
So $|p_1(A)| \leq |p_1(A_i)| -2  = 4$ and hence 
we obtain negativity:
$$-2p_1(P) - 3 \dim H^2_{\mu} \leq 8 -9 =-1$$
which  gives a contradiction.
This completes the proof of theorem $0.2$.

\vspace{3mm}

\section{Some aspects on global analysis of moduli spaces}
{\bf 6.A Deformation on  self dual curvature functionals:}
 We assume that $X$ is compact in $6.A$.
 
Let us recall the smooth function $f_{\epsilon}$ and 
  the regularization in $4.C$.
  Let 
 $A$ be an  $L^2$ ASD connection.
For 
   $a \in L^2_k(X: Ad(P) \otimes \Lambda^1) \cap (\text{ Ker } d^+_A)^{\perp}$,
   the equality: 
 $$\tilde{F}^+_A(a) = 
 d_A^+(a) +  (Q_{f_{\epsilon}}(a) \wedge Q_{f_{\epsilon}}(a))^+ =
F^+(A+a)$$ holds  for sufficiently small $\epsilon >0$.
 
 Let us  introduce a family of Hilbert spaces
which are obtained by completion of  $(\text{ Ker } d^+_A )^{\perp}\subset L^2_k$ by:
$$||a||\hat{\frak L}_k(A)_{\epsilon} = ||P_{f_{\epsilon}}(d^+_A(a))||L^2_{k-1}$$
Notice that 
 the spectra of $\Delta_A$ on $L^2(X: Ad(P) \otimes \Lambda^2_+)$
  is discrete, since $X$ is assumed to be compact.
If spectrum $\delta$ of $\Delta_A$ lie between $0 < 2 \epsilon_0 < \delta <  \epsilon_1$, 
then there is a finite dimensional vector space $H_{\epsilon_0, \epsilon_1}$
such that the isomorphism holds:
$$\hat{\frak L}_k(A)_{\epsilon_1} \oplus H_{\epsilon_0, \epsilon_1} \cong
\hat{\frak L}_k(A)_{\epsilon_0}.$$

Let $\delta_1 >0$ be the first eigenvalue, and choose $\epsilon_0  $ and $\epsilon_1$ as above.
Notice the equalities:
 $$\hat{\frak L}_k(A) _{\epsilon_0} = \hat{\frak L}_k(A) \equiv  (\text{ Ker }d^+_A)^{\perp}
  \subset L^2_k(X: Ad(P) \otimes \Lambda^1).$$
For $0< t \leq  \epsilon_1$, let us consider  the  smooth path $f_t$, and 
define:
\begin{align*}
& \tilde{F}^+_t : \hat{\frak L}_k(A) \to L^2_{k-1} (X: Ad(P) \otimes \Lambda^2_+) \\
& \tilde{F}^+_t(a) = 
\begin{cases}
d^+_A(a) + (Q_{f_t}(a) \wedge Q_{f_t}(a))^+ &  t \leq  \epsilon_1 \\
d^+_A(a) + (Q_{f_t}(a') \wedge Q_{f_t}(a'))^+ &  t \geq  \epsilon_1 \\
\end{cases}
\end{align*}
where $a=a'+a''  \in \hat{\frak L}_k(A)_{\epsilon_1} \oplus H_{\epsilon_0, \epsilon_1}$.
This is a uniformly bounded family of functionals.

\begin{lem} 
Suppose $A$ is regular. Then there is a smooth path 
$a_t \in {\frak L}_k(A)$ so that 
they satisfy the solutions 
$\tilde{F}^+_t(a_t)=0$ with $a_0=0$.

In particular there is a canonical path between $A$ and a solution to the regularized equation
$\tilde{F}^+(a)=0$.
\end{lem}
{\em Proof:}
This  follows by lemma $4.3$ and the above deformation of the ASD equations.
This completes the proof.
\vspace{3mm} \\
{\em Remark 6.1:}
Let $(S,h)$ be the Riemannian-Casson handle  homogeneously of bounded type.
Let $A$ be a smooth $L^2$ ASD connection over $S$, and consider a small neighborhood 
$U \subset   A +  L^2_k(S; Ad(P) \otimes  \Lambda^1 )$.
Then  the codimension of:
$$d^+_{A'}: {\frak L}_k(A') \to  L^2_{k-1}(S; Ad(P) \otimes \Lambda^2_+)$$
 is finite dimensional
for any $A' \in U$.

Let us denote
$codim(A') =  codim \{ 
d^+_{A'}: {\frak L}_k(A') \to L^2_{k-1}(S; Ad(P) \otimes \Lambda^2_+)\} $, and put:
$$codim(U) = \sup_{A' \in U} codim (A')< \infty.$$
If $ codim (U) = codim (A)$ holds, then 
Then  there is some $N$ such that for any $A' \in U$, there is $N$ such that
the isomorphisms 
$I_{A'} : {\frak L}_k(A) \oplus {\mathbb R}^N \cong {\frak L}_k(A')$ hold.
\vspace{3mm} \\
{\bf 6.B Yang-Mills functional:}
Let us  say that 
 $A' \in {\frak L}_k(A)$ is ASD, if 
  there is a convergent sequence $A_i' \to A' \in {\frak L}_k(A)$
 with:  
 $$A_i'  \in A + L^2_k(X: Ad(P) \otimes \Lambda^1) , \quad F^+_{A_i'}\to 0 \in L^2_{k-1}.$$
  It would be of interest for us to ask existence of   an ASD element $A' \in {\frak L}_k(A)$
  such that $A'$ is not contained in the image of 
  $L^2_k(X: Ad(P) \otimes \Lambda^1) \to {\frak L}_k(A)$.

For $A' \in A+ L^2_k(X: Ad(P) \otimes \Lambda^1) $, 
the Yang-Mills functional is given by:
$$\int_S ||F^+(A')||^2 vol.$$
Notice that this cannot be directly defined on ${\frak L}_k(A)$.

For  $A' \in {\frak L}_k(A) $, 
let  us denote by $[[a_i]]$
where $a_i  \in A + L^2_k(X: Ad(P) \otimes \Lambda^1) $ 
 converge to $A'$ in   ${\frak L}_k(A) $.
Then the Yang-Mills functional on ${\frak L}_k(A)$ is defined by:
$${\em YM}_A(A') = \inf_{[[a_i]] } \  \lim \inf_i  \ \int_S ||F^+(a_i)||^2 vol.$$

When the bundle  $E \to S$ admits  a minimal  ASD connection
and if  ${\em YM}_A(A') =0$ holds for some $A' \in {\frak L}_k(A)$,
then it would be of interest for us to ask whether $A'$ lies in the image of 
$L^2_k$ by use of sub lemma $3.1$.

\vspace{1cm}

Tsuyoshi Kato

Department of Mathematics

Faculty of Science

Kyoto University

Kyoto 606-8502
Japan

\end{document}